\documentclass[11pt]{article}

\usepackage{geometry}
\date{}

\usepackage{amsthm}
\usepackage{graphics} 
\usepackage{epsfig} 
\usepackage{times}
\usepackage{amsmath} 
\usepackage{amssymb}  
\usepackage{color}
\newtheorem{assumption}{Assumption}
\newtheorem{lemma}{Lemma}

\newtheorem{remark}{Remark}
\newtheorem{definition}{Definition}

\newtheorem{theorem}{Theorem}

\usepackage{bbm}
 \usepackage{url}
\usepackage{amssymb} 
\usepackage{amsmath}
\DeclareMathOperator*{\argmin}{argmin}

\usepackage[utf8]{inputenc} 
\usepackage[T1]{fontenc}    
\usepackage{hyperref}       
\usepackage{url}            
\usepackage{booktabs}       
\usepackage{amsfonts}       
\usepackage{nicefrac}       
\usepackage{microtype}      
\usepackage{xcolor}         
\usepackage[font=small,labelfont=bf]{caption}               
\usepackage{wrapfig}
\usepackage{float}
\usepackage{algorithm}
\usepackage{algpseudocode}
\usepackage{thmtools, thm-restate}

\usepackage{caption}
\usepackage{subcaption}
\usepackage{enumitem}

\usepackage{ifthen}
\newif\ifshortver
\usepackage{caption}

\hypersetup{
    unicode=false,          
    pdftoolbar=true,        
    pdfmenubar=true,        
    pdffitwindow=false,     
    pdfstartview={FitH},    
    pdfnewwindow=true,      
    colorlinks=true,       
    linkcolor=blue,          
    citecolor=black,        
    filecolor=magenta,      
    urlcolor=blue,           
    breaklinks=true
}

\newcommand{\mc}{\mathbb}

\usepackage{tikz}

\title{Meta-Learning Linear Quadratic Regulators:\\
A Policy Gradient MAML Approach for Model-free LQR}

\usepackage[round]{natbib}

\allowdisplaybreaks
\usepackage[toc,page,header]{appendix}
\usepackage{minitoc}

\author{Leonardo F. Toso, Donglin  Zhan, James Anderson, and Han Wang
\thanks{Leonardo F. Toso, Donglin  Zhan, James Anderson, and Han Wang are with the Department of Electrical Engineering, Columbia University in
the City of New York. Email: \texttt{\{lt2879, dz2478, james.anderson, hw2786\}@columbia.edu}.}}

\begin{document}

\doparttoc 
\faketableofcontents 


\maketitle

\begin{abstract}
We investigate the problem of learning linear quadratic regulators (LQR) in a multi-task, heterogeneous, and model-free setting. We characterize the stability and personalization guarantees of a policy gradient-based (PG) model-agnostic meta-learning (MAML) \citep{finn2017model} approach for the LQR problem under different task-heterogeneity settings. We show that our \texttt{MAML-LQR} algorithm produces a stabilizing controller close to each task-specific optimal controller up to a task-heterogeneity bias in both model-based and model-free learning scenarios. Moreover, in the \emph{model-based} setting, we show that such a controller is achieved with a linear convergence rate, which improves upon sub-linear rates from existing work. Our theoretical guarantees demonstrate that the learned controller can efficiently \emph{adapt} to unseen LQR tasks.
\end{abstract}

\section{Introduction}

One of the main successes of reinforcement learning (RL) (for example, in the context of robotics) is its ability to learn control policies that rapidly adapt to different agents and environments \citep{wang2016learning, duan2016rl,rothfuss2018promp}. This idea of learning a control policy that efficiently adapts to unseen RL tasks is referred to as meta-learning, or learning to learn. The  most popular approach is the model-agnostic meta-learning (MAML) \citep{finn2017model,finn2019online}. In the context of RL, the role of MAML is to exploit task diversity of RL tasks drawn from a common task distribution to learn a control policy in a multi-task and heterogeneous setting that is only a few policy gradient (PG) steps away from an unseen task's optimal policy.

Despite its  success in image classification and RL, more needs to be understood about the theoretical convergence guarantees of MAML for both model-based and model-free learning. This is due to the fact that, in general, the MAML objective is non-convex and requires a careful analysis depending on the considered task-setting (e.g., classification and regression \citep{fallah2020convergence,johnson2011safe,abbas2022sharp,ji2022theoretical,zhan2024data}, RL \citep{fallah2021convergence,liu2022theoretical, beck2023survey}). There is a recent body of work on multi-task/agent learning for estimation \citep{zhang2023multi,wang2023fedsysid,zhang2024sampleefficient,toso2023learning, chen2023multi} and control \citep{wang2023model,tang2023zeroth,wang2023fleet, toso2024asynchronous}, where theoretical guarantees are provided for different learning techniques in a variety of control settings. Therefore, to characterize the personalization guarantees of MAML for a baseline and well-established control setting, we consider the model-free MAML-LQR problem.

In the optimal control domain, a highly desired feature is the fast adaptation of a designed controller to unseen situations during deployment, for example, in the setting where a manufacturer  (of e.g., robots or drones) is responsible for designing optimal controllers for individual systems' objectives. Designing such controllers from scratch is sample-inefficient since it requires a large amount of trajectory data and several PG steps. Since manufacturers utilize production lines with the purpose of producing nearly identical systems, a controller should not need to be designed from scratch for every system. As such, the MAML-LQR approach exploits this similarity amongst systems within the same fabrication slot to design an LQR controller that \emph{adapts} to fresh new slots of systems. This significantly reduces the amount of trajectory data required since it relies now on a simple fine-tuning step of the learned MAML-LQR controller to suit each system's objective.

Even for a simple discrete-time control setting, provably guaranteeing that a MAML-LQR approach produces a controller that adapts to unseen LQR tasks is not an easy endeavor and requires careful handling of the task heterogeneity and the stability of the sampled tasks under the learned controller. As well-established in the literature of PG methods for the LQR problem \citep{fazel2018global, malik2019derivative, gravell2020learning, mohammadi2019global, hu2023toward}, some properties of the LQR objective (e.g., gradient dominance and local smoothness) are crucial to derive global convergence guarantees. Although tempting, we cannot simply extend these guarantees to the MAML-LQR approach since those properties of the LQR cost are no longer valid for the meta learning objective \citep{molybog2021does,musavi2023convergence}.

In contrast to \cite{molybog2021does, musavi2023convergence}, this work establishes personalization guarantees for the MAML-LQR problem in both model-based and model-free settings. In particular, \cite{molybog2021does}  only characterizes the convergence of the \emph{model-based} learning for the \emph{single-task} setting. Here, we consider the multi-task and heterogeneous setting. \cite{musavi2023convergence} establish convergence to a stationary point (i.e., local convergence analysis). It is not clear in their results how the heterogeneity across the tasks may impact the convergence of the MAML approach and how efficiently the learned controller adapts to unseen tasks. Therefore, in this work, we address these points and provide meaningful personalization guarantees that support the ability of the learned controller to adapt to unseen LQR tasks under different task-heterogeneity settings. It is also worth emphasizing that our MAML-LQR setting is different from the one considered in \cite{richards2023control}, where a control-oriented meta-learning approach is proposed to design adaptive control laws for the nonlinear feedback control problem.

\textbf{Contributions:} Toward this end, our main contributions are summarized as follows:
\begin{itemize}
\vspace{-0.3cm}
    \item  This is the first work to provide personalization guarantees for both model-based and model-free learning settings (Theorems \ref{theorem:convergence_model_based} and \ref{theorem:convergence_model_free}). Our convergence bounds characterize the distance between the learned and MAML-LQR optimal controller to each task-specific optimal controller and reveal the ultimate goal of the MAML-LQR approach, i.e., the quick adaptation to unseen tasks. In the \emph{model-based} setting, we show that the learned controller is achieved with a linear convergence rate that improves upon sub-linear rates in the existing work.
\vspace{-0.3cm}   
    \item This is the first work to establish stability (Theorem \ref{theorem:stability_model_free}) and convergence guarantees (Theorem \ref{theorem:convergence_model_free}) for the MAML-LQR approach in the model-free setting. Our convergence guarantees demonstrate that the learned controller stabilizes and is close to each task-specific optimal controller up to a task-heterogeneity bias. Furthermore, our analysis underscores the impact of different heterogeneity settings (i.e., system heterogeneity, cost heterogeneity, system, and cost heterogeneity) on the convergence of the MAML-LQR.
\end{itemize}

\vspace{-0.3cm}
\section{Model Agnostic Meta-Learning (MAML) for the LQR problem}

Consider $M$ discrete-time and linear time-invariant (LTI) dynamical systems 
\begin{align}\label{eq:LTI_system}
    x^{(i)}_{t+1} = A^{(i)}x^{(i)}_{t} + B^{(i)}u^{(i)}_t, \quad t = 0,1,2,\ldots,
\end{align}
where $A^{(i)} \in \mathbb{R}^{n_x \times n_x}$, $B^{(i)} \in \mathbb{R}^{n_x \times n_u}$, with $n_x\ge n_u$. The initial state of \eqref{eq:LTI_system} is drawn from an arbitrary distribution $\mathcal{X}_0$ that satisfies\footnote{The expectation is with respect to $x^{(i)}_0 \sim \mathcal{X}_0$.} $\mc{E}[x_0^{(i)}]=0$  and $\mc{E}[x_0^{(i)}x_0^{(i)\top}] \succ \mu I_{n_x}$ , for some $\mu > 0$, for all $i \in [M]$.\footnote{This assumption is standard in PG methods for the LQR problem \citep{fazel2018global, malik2019derivative}. It guarantees that all stationary solutions are global optima.}  The objective of the LQR problem is to design an optimal control sequence $u^{(i)}_t := - K^\star_i x^{(i)}_t$ that minimizes a quadratic cost in both states $x^{(i)}_t$ and input $u^{(i)}_t$. The optimal controllers $K^\star_i$ solve
\begin{align} \label{eq:LQR_cost}
K^\star_i :=  & \argmin_{K \in \mathcal{K}^{(i)}} \left\{ J^{(i)}(K) := \mc{E} \left[\sum_{t=0}^{\infty} x^{(i)\top}_t \left(Q^{(i)} + K^\top R^{(i)} K\right) x^{(i)}_t\right] \right\},\;\ \text{s.t}\:\ \eqref{eq:LTI_system},
\end{align}
where $Q^{(i)} \in \mathbb{S}^{n_x}_{\succ 0}$, $R^{(i)} \in \mathbb{S}^{n_u}_{\succ 0}$, and $\mathcal{K}^{(i)} := \{K~|~ \rho(A^{(i)} - B^{(i)}K) < 1\}$ denotes the set of stabilizing controllers of the $i^{\mathrm{th}}$ system, and $\rho(\cdot)$ denotes the spectral radius.
\begin{definition} The LQR task is a tuple $\mathcal{T}^{(i)} := (A^{(i)}, B^{(i)}, Q^{(i)}, R^{(i)})$ equipped with the objective of designing $K^\star_i$ that minimizes the LQR cost $J^{(i)}(K)$. 
\end{definition}

 Consider a distribution of LQR tasks denoted by $p(\mathcal{T})$ from which a collection of $M$ LQR tasks $\mathcal{T}:=\{\mathcal{T}^{(i)}\}_{i=1}^M$ are sampled. The objective of the MAML approach for the LQR problem is to design a controller $K^\star_{\text{ML}}$ based only on the tasks in $\mathcal T$ that can \emph{efficiently adapt} to any unseen LQR task originating from $p(\mathcal{T})$, i.e., we aim to find a controller that is only a few PG iterations away from any unseen task-specific optimal controller. Precisely, $K^\star_{\text{ML}}$ solve 
\begin{align} \label{eq:MAML_LQR_cost}
K^\star_{\text{ML}} :=  & \argmin_{K \in \bar{\mathcal{K}}} \left\{ J_{\text{ML}}(K) := \frac{1}{M}\sum_{i=1}^M J^{(i)}\underbrace{\left(K - \eta_l \nabla J^{(i)}(K)\right)}_{\mathrm{single\text{-}step~PG}} \right\}, \;\ \text{s.t}\:\ \eqref{eq:LTI_system} \;\ \forall i \in [M],
\end{align}
where $\bar{\mathcal{K}}:= \cap_{i\in [M]} \mathcal{K}^{(i)}$ is the MAML stabilizing set and $\eta_l$ denotes some positive step-size. To solve \eqref{eq:MAML_LQR_cost}, we exploit a PG-based approach where the update rule is described as follows:
\begin{align}\label{eq:MAML_LQR_updating_rule}
    K \leftarrow K -\eta \nabla J_{\text{ML}}(K), \text{ where } \nabla J_{\text{ML}}(K) := \frac{1}{M}\sum_{i=1}^{M} H^{(i)}(K) \nabla J^{(i)}(K - \eta_l \nabla J^{(i)}(K)),
\end{align}
with $H^{(i)}(K) := I_{n_u} -\eta_l \nabla^2 J^{(i)}(K)$, and $\eta$ being some positive (possibly time-varying) step-size. Next, we define the task-specific and MAML stabilizing sub-level sets.  

\begin{definition} (Stabilizing sub-level set) \label{def:stabilizing_set} The task-specific and MAML stabilizing sub-level sets are defined as follows:

\begin{itemize}
    \item Given a task $\mathcal{T}^{(i)}$, the task-specific sub-level set $\mathcal{S}^{(i)} \subseteq \mathcal{K}^{(i)}$ is
    \begin{align*}
        \mathcal{S}^{(i)}:= \left\{K\; | \; J^{(i)}(K) - J^{(i)}(K^\star_i) \leq \gamma_i \Delta^{(i)}_0\right\}, \text{ with } \Delta^{(i)}_0 = J^{(i)}(K_0) - J^{(i)}(K^\star_i).  
    \end{align*}
    where $K_0$ denotes an initial controller and $\gamma_i$ being any positive constant.

    \item The MAML-LQR stabilizing sub-level set ${\mathcal{S}}_{\text{ML}} \subseteq \bar{\mathcal{K}}$ is defined as the intersection between each task-specific stabilizing sub-level set, i.e., ${\mathcal{S}}_{\text{ML}} := \cap_{i\in [M]} \mathcal{S}^{(i)}$. 
\end{itemize}
\end{definition}
\begin{remark}
 Observe that, if $K\in {\mathcal{S}}_{\text{ML}}$, i.e., $K$ stabilizes all LQR tasks in $\mathcal{T}$, one may select a step-size $\eta_l$, such that $\bar{K} = K - \eta_l \nabla J^{(i)}(K)$ also stabilizes all the LQR tasks in $\mathcal{T}$, i.e., $\bar{K}\in {\mathcal{S}}_{\text{ML}}$. We prove this fact and provide the condition on $\eta_l$ to satisfy it in the stability analysis in Theorem \ref{theorem:stability_model_based}. 
\end{remark}

\begin{assumption} \label{assumption:initial_stabilizing_K0} We have access to an initial stabilizing controller $K_0 \in \mathcal{S}_{\text{ML}}$. 
\end{assumption}

\begin{remark}
The above assumption is standard in PG methods for the LQR problem \citep{fazel2018global, gravell2020learning,wang2023model,toso2023oracle}.  If the initial controller $K_0$ fails to stabilize \eqref{eq:LTI_system}, $\forall i\in [M]$, the MAML-LQR update in \eqref{eq:MAML_LQR_updating_rule} cannot produce a stabilizing controller, since ${\nabla} J^{(i)}(K_0)$, ${\nabla}^2 J^{(i)}(K_0)$ are both undefined for the corresponding unstabilized tasks. \cite{perdomo2021stabilizing} and \cite{ozaslan2022computing} detail how to find an initial stabilizing controller for the single LQR instance. Moreover, it is worth emphasizing that although $K_0$ stabilizes \eqref{eq:LTI_system}, $\forall i\in [M]$, it may provide a sub-optimal performance, i.e., $J^{(i)}(K_0) \geq J^{(i)}(K^\star_i)$. 
\end{remark}

\vspace{-0.35cm}

As well-established in the literature of PG-LQR, $J^{(i)}(K)$ is, in general, non-convex with respect to $K$. However, by leveraging some properties of the LQR cost (e.g., gradient domination and local smoothness), \cite{fazel2018global} provide global convergence guarantees of PG methods for both model-based and model-free LQR settings. Although tempting, these properties of the LQR cost cannot simply be extended to the MAML-LQR objective when dealing with task heterogeneity as discussed in \cite{molybog2021does}.

In the sequel, we proceed as follows: We first provide conditions on the problem parameters to ensure that given any stabilizing controller $K \in \mathcal{S}_{\text{ML}}$, $K - \eta \nabla J_{\text{ML}}(K)$ is also MAML stabilizing, i.e.,  $K - \eta \nabla J_{\text{ML}}(K) \in \mathcal{S}_{\text{ML}}$. In contrast to \cite{musavi2023convergence}, which guarantees that a model-based MAML-LQR approach finds a stationary solution, we derive global convergence bounds for model-based and model-free learning. In particular, our convergence guarantees underscore the impact of different task-heterogeneity settings on the closeness of the learned controller and each task-specific optimal controller and demonstrate its adaptation to unseen tasks.
\vspace{-0.3cm}
\subsection{Model-based LQR}

In the model-based LQR setting, we assume to have access to the tuple $\mathcal{T}^{(i)} = (A^{(i)}, B^{(i)}, Q^{(i)}, R^{(i)})$. With the ground-truth model in hand, we have closed-form expressions to compute both gradient $\nabla J^{(i)}(K)$ and Hessian $\nabla^2 J^{(i)}(K)$ of the LQR cost:

\begin{itemize}
\vspace{-0.2cm}
    \item \textbf{Gradient of the LQR cost $\nabla J^{(i)}(K)$} \citep{fazel2018global}: Given $\mathcal{T}^{(i)}$ and a stabilizing controller $K \in \mathcal{S}_{\text{ML}}$, the gradient is given by $\nabla J^{(i)}(K):=2E^{(i)}_K\Sigma^{(i)}_{K}$, where 
   \begin{align*}
    E^{(i)}_K & := R^{(i)}K - B^{(i)\top} P^{(i)}_K(A^{(i)} - B^{(i)} K), \quad \text{and} \quad \Sigma^{(i)}_K := \mc{E}\sum_{t=0}^\infty x^{(i)}_t x^{(i)\top}_t,
    \end{align*}    
   with $x^{(i)}_t$ subject to the system dynamics in \eqref{eq:LTI_system}, and $P^{(i)}_K \in \mathbb{S}^{n_x}_{\succ 0}$ denoting the solution of the Lyapunov equation $P^{(i)}_{K} := Q^{(i)} +K^{\top}R^{(i)}K + (A^{(i)}-B^{(i)}K)^{\top}P^{(i)}_{K}(A^{(i)}-B^{(i)}K)$. 
\vspace{-0.2cm}    
    \item \textbf{Hessian of the LQR cost $\nabla^2 J^{(i)}(K)$} \citep{bu2019lqr}: Given $\mathcal{T}^{(i)}$ and a stabilizing controller $K \in \mathcal{S}_{\text{ML}}$, the Hessian operator at $K$ acting on some $X \in \mathbb{R}^{n_u\times n_x}$, is given by
\begin{align*}
\nabla^2 J^{(i)}(K)[X] :=  2\left(R^{(i)}+B^{(i)\top} P^{(i)}_K B^{(i)}\right) X \Sigma^{(i)}_K  -4 B^{(i)\top} \tilde{P}_K^{(i)}[X](A^{(i)}-B^{(i)} K) \Sigma^{(i)}_K,
\end{align*}
with $\tilde{P}^{(i)}_K[X]:=(A^{(i)}-B^{(i)} K)^\top \tilde{P}_K^{(i)}[X](A^{(i)}-B^{(i)} K)+X^\top E^{(i)}_K+  E^{(i)\top}_K X$.

\end{itemize}

Hence, we can exploit these closed-form expressions in order to update the controller (see step~\ref{eq:MAML_LQR_updating_rule} in \texttt{MAMl-LQR}, Algorithm~\ref{alg:MAML_LQR}) for the model-based MAML-LQR. In step 4, our algorithm computes a one-step inner gradient descent iteration on $K_n$ and $H^{(i)}(K_n) = I_{n_u} - \eta \nabla^2 J^{(i)}(K_n)$, for each task $i \in [M]$ and iteration $n$. These quantities are then used to update the controller $K_{n+1}$ in step 6. By repeating these steps for $N$ iterations, Algorithm \ref{alg:MAML_LQR} returns $K_N$. We further prove that $K_N$ is close to each task-specific optimal controller $K^{\star}_i$, which in turn is proved to be close to the MAML-LQR optimal controller $K^{\star}_{\text{ML}}$. To demonstrate this, we  revisit some properties of the LQR cost function.

\begin{lemma}\label{lemma:uniform_bounds} (Uniform bounds) Given $\mathcal{T}^{(i)}$ and a stabilizing controller $K \in \mathcal{S}_{\text{ML}}$, the  gradient $\nabla J^{(i)}(K)$, Hessian $\nabla^2 J^{(i)}(K)$, and controller $K$ are bounded as follows:
$$
\|\nabla J^{(i)}(K)\|_F \leq h_G(K), \quad \|\nabla^2 J^{(i)}(K)\|_F \leq h_H(K), \quad \text {and}\quad \|K\|_F \leq h_c(K),
$$
where $h_G(K)$, $h_H(K)$, and $h_c(K)$ are functions of the problem parameters.
\end{lemma}

\begin{lemma}\label{Lemma:lipschitz} (Local smoothness) Given  $\mathcal{T}^{(i)}$ and two stabilizing controllers $K, K^{\prime} \in \mathcal{S}_{\text{ML}}$ such that $\|\Delta\| := \|K^{\prime} -K\|  \leq h_{\Delta}(K)\ <\infty$.  The LQR cost, gradient and Hessian satisfy:
\begin{align*}
&\left|J^{(i)}\left(K^{\prime}\right)-J^{(i)}(K)\right| \leq h_{\text {cost}}(K) J^{(i)}(K) \|\Delta\|_F, \\
&\left\|\nabla J^{(i)}\left(K^{\prime}\right)-\nabla J^{(i)}(K)\right\|_F \leq h_{\text {grad}}(K)\|\Delta\|_F, \\
&\left\|\nabla^2 J^{(i)}\left(K^{\prime}\right)-\nabla^2 J^{(i)}(K)\right\|_F \leq h_{\text {hess }}(K)\|\Delta\|_F
\end{align*}
 where $h_{\Delta}(K)$, $h_{\text {cost}}(K)$, $h_{\text {hess }}(K)$ and $h_{\text {grad}}(K)$ are functions of the problem parameters.
\end{lemma}

\begin{lemma}\label{lemma:gradient_domination} (Gradient Domination) Given $\mathcal{T}^{(i)}$ and a stabilizing controller $K \in \mathcal{S}_{\text{ML}}$. Let $K^\star_i$ be the optimal controller of task $\mathcal{T}^{(i)}$. Then, it holds that
\begin{align*}
     J^{(i)}(K)-J^{(i)}\left(K^\star_i\right) \leq \frac{1}{\lambda_i}\|\nabla J^{(i)}(K)\|_F^2
\end{align*}
where $\lambda_i := 4\mu^2 \sigma_{\min }(R^{(i)})/\left\|\Sigma_{K^\star_i}\right\|$. 
\end{lemma}

The uniform bounds of $\|\nabla J^{(i)}(K)\|_F$ and $\|K\|_F$, and the gradient domination property are proved in \citep{fazel2018global, wang2023fedsysid}. Moreover, the uniform bound of $\|\nabla^2 J^{(i)}(K)\|_F$ can be found in \citep[Lemma 7.9]{bu2019lqr}. In addition, the proofs for the local smoothness of the cost and gradient are detailed in \citep[Appendix F]{wang2023model}, whereas the local smoothness of the Hessian is proved in \citep[Appendix B]{musavi2023convergence}. The explicit expressions of $h_G(K)$, $h_c(K)$, $h_H(K)$, $h_{\Delta}(K)$, $h_{\text {cost}}(K)$, and $h_{\text {grad}}(K)$ are revisited in Appendix \ref{appendix:polynomials}. Throughout the paper, we use $\bar{h}:= \sup_{K \in \mathcal{S}_{\text{ML}}} {h}(K)$ and $\underline{h}:= \inf_{K \in \mathcal{S}_{\text{ML}}} {h}(K)$ to denote the supremum and infimum of some positive polynomial $h(K)$ over the set of stabilizing controllers ${\mathcal{S}_{\text{ML}}}$.

\begin{algorithm}
\caption{ \texttt{MAML-LQR}: Model-Agnostic Meta-Learning for LQR tasks \textbf{(Model-based)}} 
\label{alg:MAML_LQR}
\begin{algorithmic}[1]
\State \textbf{Input:} initial stabilizing controller $K_0$, inner and outer step-sizes $\eta_l$, $\eta$
\State \textbf{for} $n=0, \ldots, N-1$ \textbf{ do} 
\State \quad  \textbf{for} each task $i \in [M]$ in $\mathcal{T}$ \textbf{compute}
\State \quad \quad $\bar{K}^{(i)}_n = K_n - \eta_l \nabla J^{(i)}(K_n)$, \text{and} $H^{(i)}(K_n) = I_{n_u} - \eta_l\nabla^2J^{(i)}(K_n)$
\State \quad  \textbf{end for}
\State \quad \quad     $K_{n+1} = K_{n} - \frac{\eta}{M}\sum_{i=1}^M H^{(i)}(K_n) \nabla J^{(i)}(\bar{K}^{(i)}_n)$
\State\textbf{end for}
\State \textbf{Output:} $K_{N}$
\end{algorithmic}
\end{algorithm}

\vspace{-0.4cm}
\subsection{Task Heterogeneity}

In contrast to \cite{musavi2023convergence}, we consider the MAML-LQR problem in a variety of heterogeneity settings. Our goal is to determine precisely how task heterogeneity impacts convergence of \texttt{MAML-LQR} -- in model-based setting (Algorithm~\ref{alg:MAML_LQR}) and the \emph{model-free setting} (described later). We consider a task-heterogeneity setting characterized by the combination of system and cost heterogeneity. That is, we assume that there exist positive scalars $\epsilon_1$, $\epsilon_2$, $\epsilon_3$ and $\epsilon_4$, such that $$\underset{i\neq j}{\max} \lVert  A^{(i)} -A^{(j)}\rVert \leq \epsilon_1, 
\underset{i\neq j}{\max}\lVert B^{(i)} -B^{(j)} \rVert\leq \epsilon_2, \underset{i\neq j}{\max} \lVert  Q^{(i)} - Q^{(j)}\rVert \leq \epsilon_3,  
\underset{i\neq j}{\max}\lVert R^{(i)} -R^{(j)} \rVert\leq \epsilon_4.$$

Observe that this setting spans three different types of task heterogeneity: 1) system heterogeneity, with $\epsilon_3=\epsilon_4=0$, i.e., $Q^{(i)} = Q$, $R^{(i)} = R$, $\forall i \in [M]$. 2) cost heterogeneity, with $\epsilon_1=\epsilon_2=0$, i.e., $A^{(i)} = A$, $B^{(i)} = B$, $\forall i \in [M]$. 3) system and cost heterogeneity, where $\epsilon_1$, $\epsilon_2$, $\epsilon_3$ and $\epsilon_4$ are non-zero. Next, we bound the norm of the gradient difference between two distinct tasks.

\vspace{-0.2cm}
\begin{lemma}\label{lemma:gradient_heterogeneity} (Gradient heterogeneity) For any two distinct LQR tasks $\mathcal{T}^{(i)}$ and $\mathcal{T}^{(j)}$, and stabilizing controller $K \in \mathcal{S}_{\text{ML}}$. It holds that, 
\begin{align}\label{eq:gradient_heterogeneity}
    \|\nabla J^{(i)}(K) - \nabla J^{(i)}(K)\| \leq {f}_z(\bar{\epsilon}) := \epsilon_1 {h}^1_{\text{het}}(K) + \epsilon_2 {h}^2_{\text{het}}(K) + \epsilon_3 {h}^3_{\text{het}}(K) + \epsilon_4 {h}^4_{\text{het}}(K),
\end{align}
for any $i\neq j \in [M]$, where $z$ $\in \{1,2,3\}$,\footnote{ $z=1$ refers to the system heterogeneity, $z=2$ to cost heterogeneity and $z=3$ to system and cost heterogeneity.} and $\Bar{\epsilon} = \{\epsilon_1, \epsilon_2, \epsilon_3, \epsilon_4\}$, where ${h}^1_{\text{het}}(K), {h}^2_{\text{het}}(K), {h}^3_{\text{het}}(K)$, and ${h}^4_{\text{het}}(K)$ are positive polynomials that depend on the problem parameters. 
\end{lemma}

The proof and explicit expressions of ${h}^1_{\text{het}}(K), {h}^2_{\text{het}}(K), {h}^3_{\text{het}}(K)$, and ${h}^4_{\text{het}}(K)$, are detailed in Appendix \ref{appendix:proof_gradient_het}. We observe that as long as the heterogeneity level $\epsilon_1, \epsilon_2, \epsilon_3, \epsilon_4$ is small, the gradient descent direction of task $\mathcal{T}^{(i)}$ \eqref{eq:MAML_LQR_updating_rule} is close to the one of task $\mathcal{T}^{(j)}$. Moreover, by combining \eqref{eq:gradient_heterogeneity} along with the uniform bound of the Hessian (Lemma \ref{lemma:uniform_bounds}), we observe that $\nabla J_{\text{ML}}(K)$ is also close to $\nabla J^{(i)}(K)$. This is a crucial step we use in our stability and convergence analysis for both model-based and model-free settings. 
\vspace{-0.4cm}
\subsection{Model-free LQR}

We now consider the setting where the tuple $\mathcal{T}^{(i)} = (A^{(i)}, B^{(i)}, Q^{(i)}, R^{(i)})$ is unknown. Therefore, computing the gradient and Hessian through closed-form expressions is no longer possible. This forces us to resort to  methods that approximate such quantities. Following numerous work in the literature of model-free PG-LQR \citep{fazel2018global,malik2019derivative,gravell2020learning,mohammadi2019global,wang2023model}, we focus on zeroth-order methods to estimate the gradient and Hessian of the LQR cost. In particular, we consider a two-point estimation scheme since it has a lower estimation variance compared to its one-point counterpart \citep{malik2019derivative}. 

Zeroth-order methods with two-point estimation solely rely on querying cost values at symmetric perturbed controllers to construct a biased estimation of both the gradient and Hessian. In particular, zeroth-order estimation is a Gaussian smoothing approach \citep{nesterov2017random} based on Stein's identity \citep{stein1972bound} that relates gradient and Hessian to cost queries.

\begin{algorithm}
\caption{\texttt{ZO2P}: Zeroth-order with two-point estimation} 
\label{alg:ZO2P}
\begin{algorithmic}[1]
\State \textbf{Input:} controller $K$, number of samples $m$ and smoothing radius $r$,
\State \textbf{for} $l=1, \ldots, m$ \textbf{ do } \
\State \quad Sample controllers $K^1_l=K+U_l$ and $K^2_l=K-U_l$,  where $U_l$ is drawn uniformly at random over matrices with Frobenius norm $r$. A cost oracle provides $J(K^1)$, $J(K^2)$ and $J(K)$
\State \textbf{end for}\
\State \quad $\widehat{\nabla}J(K)= \frac{n_x n_u}{2 r^2 m} \sum_{l=1}^{m}  (J(K^1_l) - J(K^2_l)) U_l$
\State \quad $\widehat{\nabla}^2J(K)=\frac{n^2_u}{r^2 m}  \sum_{l=1}^{m} (J(K^1_l) - J(K)) (U_lU_l^\top - I_{n_u})$
\State \textbf{Return} $\widehat{\nabla}J(K)$, $\widehat{\nabla}^2J(K)$
\end{algorithmic}
\end{algorithm}

Algorithm \ref{alg:ZO2P} describes the zeroth-order estimation of the gradient and Hessian of the LQR cost. First, a pair of symmetric perturbations  to the controller $K$ are sampled according to $K^1 = K + U$, $K^2 = K - U$ where $U \sim \mathbb{S}_r$. Here $\mathbb{S}_r$ denotes a distribution of $n_u \times n_x$ real matrices with $\|U\|_F = r$, where $r$ is the smoothing radius. The gradient is estimated via the first-order Gaussian Stein’s identity, $\mc{E}\left[\nabla J^{(i)}(K)\right] = \mc{E}\left[\frac{n_xn_u}{2r^2}(J^{(i)}(K^1) - J^{(i)}(K^2))U\right]$ \citep{mohammadi2020linear}, and the Hessian with its second-order counterpart, $\mc{E}\left[\nabla^2 J^{(i)}(K)\right] = \mc{E}\left[\frac{n^2_u}{r^2}(J^{(i)}(K^1) - J^{(i)}(K))(UU^\top - I_{n_u})\right]$, \citep{balasubramanian2022zeroth}, both with $m$ samples.

\vspace{-0.22cm}
\begin{remark} (Cost oracle) For simplicity, we assume that the true cost is provided by an oracle, as in \cite{malik2019derivative,toso2023oracle}. We emphasize that our work can be readily extended to the setting where only a finite-horizon approximation of the cost is available. That is the case since any finite-horizon approximation of the true cost is upper-bounded by its true value, with the approximation error controlled by the horizon length  \citep[Appendix B]{gravell2020learning}. 
\end{remark}

With the gradient and Hessian zeroth-order estimators, Algorithm \ref{alg:MAML_model_free} follows the same structure as Algorithm  \ref{alg:MAML_LQR}. The key differences are steps 4 and 7 where the gradient and Hessian computations are replaced by a zeroth-order estimation. Despite the estimation error in the zeroth-order method, we show that this error can be controlled so it does not impact the convergence of our \texttt{MAML-LQR} approach for both model-based and model-free settings. 

\begin{algorithm}
\caption{\texttt{MAML-LQR}: Model-Agnostic Meta-Learning for LQR tasks \textbf{(Model-free)}} 
\label{alg:MAML_model_free}
\begin{algorithmic}[1]
\State \textbf{Input:} initial stabilizing controller $K_0$, inner and outer step-sizes $\eta_l$, $\eta$, smoothing radius $r$, number of samples $m$.
\State \textbf{for} $n=0, \ldots, N-1$ \textbf{ do} 
\State \quad  \textbf{for} each task $i \in [M]$ in $\mathcal{T}$ \textbf{compute}
\State \quad \quad $\left[\widehat{\nabla} J^{(i)}(K_n), \widehat{\nabla}^2 J^{(i)}(K_n)\right] = \texttt{ZO2P}(K_n,m,r)$
\State \quad \quad $\widehat{K}^{(i)}_n = K_n - \eta_l \widehat{\nabla} J^{(i)}(K_n)$, \text{and} $\widehat{H}^{(i)}(K_n) = I_{n_u} - \eta_l\widehat{\nabla}^2J^{(i)}(K_n)$
\State \quad  \textbf{end for}
\State \quad \quad $\widehat{\nabla} J^{(i)}(\widehat{K}^{(i)}_n) = \texttt{ZO2P}(\widehat{K}^{(i)}_n,m,r)$,      $K_{n+1} = K_{n} - \frac{\eta}{M}\sum_{i=1}^M \widehat{H}^{(i)}(K_n) \widehat{\nabla} J^{(i)}(\widehat{K}^{(i)}_n)$
\State\textbf{end for}
\State \textbf{Output:} $K_{N}$
\end{algorithmic}
\end{algorithm}

\section{Theoretical Guarantees}\label{sec:theoretical_guarantees}

We now provide the theoretical guarantees for the stability and convergence of the MAML-LQR for both model-based and model-free settings, i.e., Algorithms \ref{alg:MAML_LQR} and \ref{alg:MAML_model_free}.

\subsection{Stability Analysis}

The objective of the stability analysis is to provide the conditions on the step-sizes $\eta_l$, $\eta$, heterogeneity $\bar{f}_z(\bar{\epsilon})$, and zeroth-order estimation parameters $m$ and $r$, such that for every iteration of Algorithms \ref{alg:MAML_LQR} and \ref{alg:MAML_model_free}, the currently obtained controller is MAML stabilizing, i.e., $K_n \in \mathcal{S}_{\text{ML}}$, $\forall n$.

\begin{theorem} \label{theorem:stability_model_based} (Model-based) Given an initial stabilizing controller $K_0 \in \mathcal{S}_{\text{ML}}$, suppose that the step-sizes and heterogeneity satisfy $\eta_l \leq \min\left\{ \frac{n_u}{\sqrt{2}h_H}, \frac{1}{\sqrt{2}\bar{h}_{\text{grad}}}, \frac{1}{\sqrt{12(12\bar{h}^2_{\text{grad}}n^2_u +\bar{h}_{H}^2 )}}\right\}$, $\eta \leq \frac{1}{4\bar{h}_{\text{grad}}}$ and $\bar{f}_z(\bar{\epsilon}) \leq  \sqrt{\min_{i} \frac{\lambda_i \Delta ^{(i)}_0}{288n^3_u}}$, respectively. Then, $\bar{K}^{(i)}_n, K_n \in \mathcal{S}_{\text{ML}}$, for every iteration of Algorithm \ref{alg:MAML_LQR}. 
    
\end{theorem}

\begin{theorem} \label{theorem:stability_model_free} (Model-free) Given an initial stabilizing controller $K_0 \in \mathcal{S}_{\text{ML}}$ and scalar $\delta \in (0,1)$, suppose that the step-sizes satisfy $\eta_l \leq \min \left \{ \frac{1}{\bar{h}_G}, \frac{n_u}{\bar{h}_H}, \frac{1}{\bar{h}_{\text{grad}}},  \frac{1}{\sqrt{20(12\bar{h}^2_{\text{grad}}n^2_u +\bar{h}_{H}^2 )}},\frac{1}{2} \right\}$, $ \eta \leq \frac{1}{8\bar{h}_{\text{grad}}}$. In addition, the heterogeneity, smoothing radius and number of samples satisfy $\bar{f}_z(\bar{\epsilon}) \leq \sqrt{\min_i \frac{\lambda_i\Delta^{(i)}_0}{480n^3_u}}$, $r \leq \min \left\{\underline{h}^1_r\left(\frac{\sqrt{{\psi}^{(i)}}}{2}\right), \underline{h}^2_r\left(\frac{\sqrt{{\psi}^{(i)}}}{2}\right) \right\}$, and $m \geq \max\left\{\bar{h}^1_m\left(\frac{\sqrt{{\psi}^{(i)}}}{2}, \delta\right),\bar{h}^2_m\left(\frac{\sqrt{{\psi}^{(i)}}}{2}, \delta\right) \right\}$\footnote{The expressions of the positive polynomials $\underline{h}^1_r(\cdot),\underline{h}^2_r(\cdot), \bar{h}^1_m(\cdot)$ and $\bar{h}^2_m(\cdot)$ are deferred to Appendix \ref{appendix:polynomials}.}, with ${\psi}^{(i)} := \frac{\lambda_i \Delta^{(i)}_0}{1296}$. Then, with probability, $1-\delta$, $\bar{K}^{(i)}_n, K_n \in \mathcal{S}_{\text{ML}}$, for every iteration of Algorithm \ref{alg:MAML_model_free}.
\end{theorem}

The proof of Theorems \ref{theorem:stability_model_based} and \ref{theorem:stability_model_free} are deferred to Sections \ref{appendix:stability_model_based} and \ref{appendix:stability_model_free} of the appendix. The proof strategy follows from an induction argument where the base case is the first iteration. We combine the local smoothness of each task-specific LQR cost (Lemma \ref{Lemma:lipschitz}) along with the gradient heterogeneity bound (Lemma \ref{lemma:gradient_heterogeneity}) and the definition of the MAML-LQR stabilizing sub-level set to show that $J^{(i)}(K_1) \leq J^{(i)}(K_0)$ for any $i \in [M]$. These results provide the conditions for which the learned controller $K_N$ is MAML stabilizing. This is essential to guarantee that the learned controller $K_N$ in Algorithms \ref{alg:MAML_LQR} and \ref{alg:MAML_model_free} can be promptly utilized to stabilize an unseen LQR task drawn from $p(\mathcal{T})$.

\vspace{-0.3cm}
\subsection{Convergence Analysis}

We now provide the conditions on the step-sizes $\eta_l, \eta$ and zeroth-order estimation parameters $m$, and $r$, such that we can ensure that the learned MAML-LQR controller $K_N$ is sufficiently close to each task-specific optimal controller $K^\star_i$ and to the optimal MAML controller $K^\star_{\text{ML}}$. For this purpose, we study the closeness of $K_N$ and $K^\star_i$ by bounding $J^{(i)}(K_N) - J^{(i)}(K^\star_i)$ and the closeness of $K^\star_{\text{ML}}$ with $J^{(i)}(K^\star_{\text{ML}}) - J^{(i)}(K^\star_i)$. 

\begin{theorem} \label{theorem:convergence_model_based} (Model-based) Given an initial stabilizing controller $K_0 \in \mathcal{S}_{\text{ML}}$, suppose that the step-sizes and number of iterations satisfy $\eta_l \leq \min \left \{\frac{n_u}{\sqrt{2}\bar{h}_H}, \frac{1}{\sqrt{2}\bar{h}_{\text{grad}}},  \frac{1}{\sqrt{12(12\bar{h}^2_{\text{grad}}n^2_u +{h}_{H}^2 )}}     \right\}$, and $ \eta \leq \frac{1}{4\bar{h}_{\text{grad}}}$,  $N \geq \frac{8}{\eta \lambda_i}\log\left(\frac{\Delta^{(i)}_0}{\epsilon^\prime}\right)$, respectively, for some small tolerance $\epsilon^\prime \in (0,1)$. Then, it holds that, 
\begin{align}
    &J^{(i)}(K_{N}) - J^{(i)}(K^\star_i)  \leq \epsilon^\prime + \frac{144n^3_u\bar{f}^2_z(\bar{\epsilon})}{\lambda_i}, \label{eq:convergence_model_based_1}\\
    &J^{(i)}(K^\star_{\text{ML}}) - J^{(i)}(K_i^\star) \leq \frac{96\bar{J}_{\max}n^3_uf^2_z(\bar{\epsilon})}{\mu^2 \min_i\sigma_{\min }(R^{(i)})\min_i\sigma_{\min }(Q^{(i)})} \label{eq:convergence_model_based_2}, 
\end{align}
with $\bar{J}_{\max} := \max_i J^{(i)}(K_0)$. 
\end{theorem}

\begin{theorem}\label{theorem:convergence_model_free} (Model-free) Given an initial stabilizing controller $K_0 \in \mathcal{S}_{\text{ML}}$ and scalar $\delta \in (0,1)$, suppose that the step-sizes satisfy $\eta_l \leq \min \left \{ \frac{1}{\bar{h}_G}, \frac{n_u}{\bar{h}_H}, \frac{1}{\bar{h}_{\text{grad}}},  \frac{1}{\sqrt{20(12\bar{h}^2_{\text{grad}}n^2_u +{h}_{H}^2 )}},\frac{1}{2} \right\}$ and $ \eta \leq \frac{1}{8\bar{h}_{\text{grad}}}$, and the smoothing radius satisfies  $r \leq \min \left\{\underline{h}^1_r\left(\frac{\epsilon^\prime \lambda_i}{1296}\right), \underline{h}^2_r\left(\frac{\epsilon^\prime \lambda_i}{1296}\right) \right\}$. Moreover, suppose that the  number of samples is selected according to $ m \geq \max\left\{\bar{h}^1_m\left(\frac{\epsilon^\prime \lambda_i}{1296}, \delta\right),\bar{h}^2_m\left(\frac{\epsilon^\prime \lambda_i}{1296}, \delta\right) \right\}$, and the number of iterations satisfies $N \geq \frac{8}{\eta \lambda_i}\log\left(\frac{2\Delta^{(i)}_0}{\epsilon^\prime}\right)$, for some small tolerance $\epsilon^\prime \in (0,1)$. Then, with probability $1-\delta$, it holds that, 
\begin{align}\label{eq:convergence_model_free}
    &J^{(i)}(K_{N}) - J^{(i)}(K^\star_i)  \leq \epsilon^\prime + \frac{240n^3_u\bar{f}^2_z(\bar{\epsilon})}{\lambda_i}, \text{ and } \eqref{eq:convergence_model_based_2}.
\end{align}
\end{theorem}

The proofs of Theorems \ref{theorem:convergence_model_based} and \ref{theorem:convergence_model_free} are detailed in Sections \ref{appendix:convergence_model_based} and \ref{appendix:convergence_model_free} of the appendix. The proof strategy follows from the local smoothness of the LQR cost (Lemma \ref{Lemma:lipschitz}),  gradient domination (Lemma \ref{lemma:gradient_domination}), and the gradient heterogeneity bound (Lemma \ref{lemma:gradient_heterogeneity}). The model-free setting also involves controlling the estimation error through a matrix Bernstein-type of inequality \citep{tropp2012user,gravell2020learning}.

These results characterize the convergence of the MAML-LQR for both model-based and model-free settings. We emphasize that both Algorithms \ref{alg:MAML_LQR} and \ref{alg:MAML_model_free} produce a controller $K_N$ that is provably close to each task-specific optimal controller $K^\star_i$ up to a heterogeneity bias. This indicates that under a low heterogeneity regime, $K_N$ will serve as a good initialization for any unseen task that is also drawn from $p(\mathcal{T})$ (i.e., an unseen task that satisfies the same task-heterogeneity level as the ones used in the MAML-LQR learning process). Moreover, in contrast to \cite{musavi2023convergence}, our convergence bounds \eqref{eq:convergence_model_based_1}, \eqref{eq:convergence_model_based_2} and \eqref{eq:convergence_model_free} emphasize the impact of task heterogeneity on the convergence of Algorithms \ref{alg:MAML_LQR} and \ref{alg:MAML_model_free}. In a low heterogeneity regime, where $K_N$ and $K^\star_i$, and $K^\star_{\text{ML}}$ and $K^\star_i$ are close, one may conclude that $K_N$ and $K^\star_{\text{ML}}$ are also sufficiently close. We also emphasize that, in the model-based setting, the learned controller $K_N$ is achieved with a linear convergence rate on the iteration count, which improves upon the sub-linear rate in \cite{musavi2023convergence}.

\vspace{-0.3cm}
\section{Experimental Results}

Numerical results\footnote{Code can be downloaded from: \url{https://github.com/jd-anderson/MAML-LQR}.} are now provided to illustrate and assess the convergence and personalization of the model-free MAML-LQR approach. In particular, we show that initializing from the learned MAML-LQR controller (i.e., $K_{0,\text{PG}} = K_N$) enables a model-free PG-LQR approach \citep[Section 4.2]{fazel2018global} to be close to the task-specific optimal controller within just a few PG iterations, for an unseen task. To illustrate this,  we consider an unstable modification of the Boeing system from \cite{honglecnotes} as the nominal LQR task. The nominal task is then used to generate multiple tasks. The technical details on the experimental setup are deferred to Appendix \ref{appendix:experimental_setup}.

\begin{figure}
  \centering
  \begin{minipage}{.32\textwidth}
    \centering
    \includegraphics[width=0.95\textwidth]{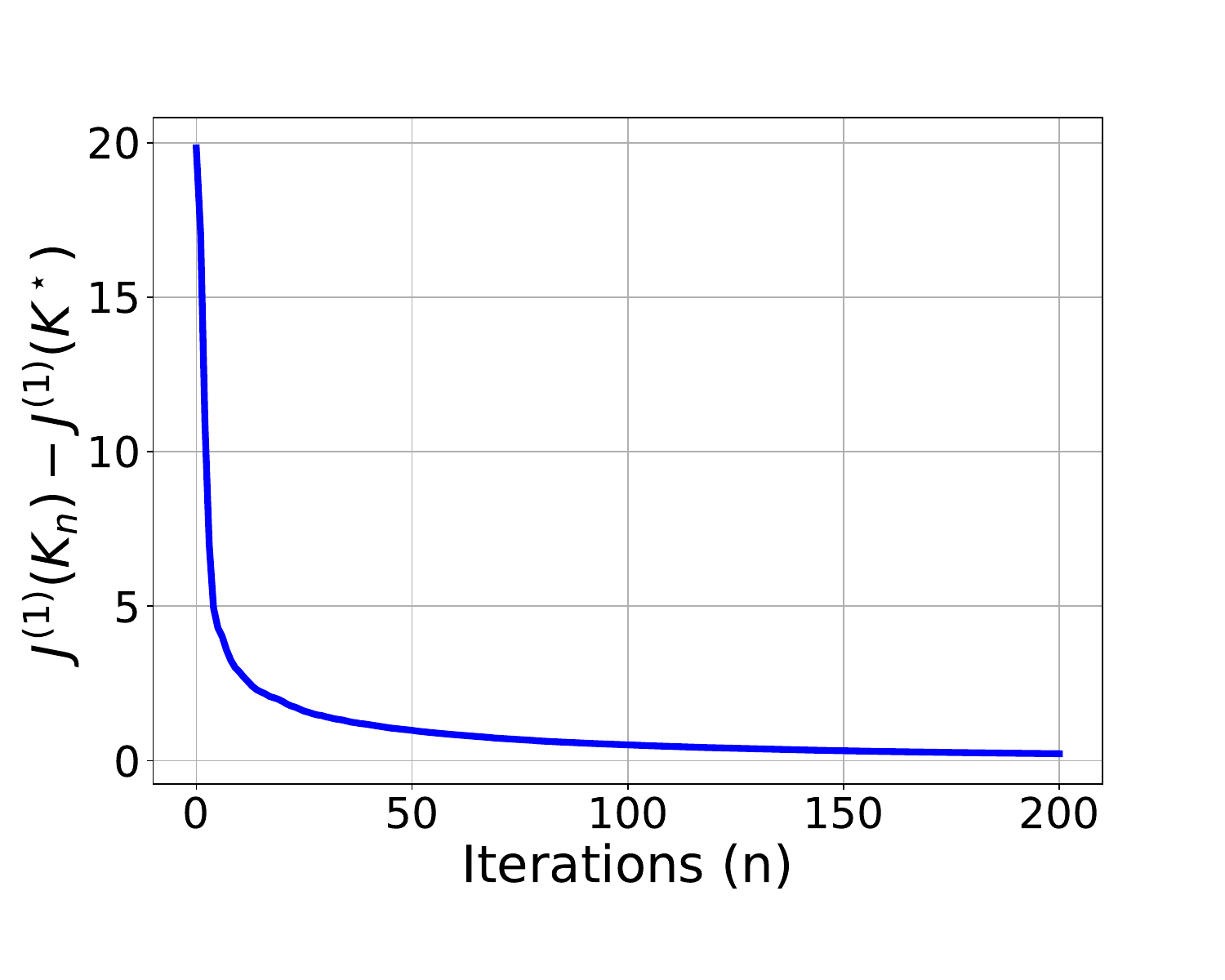}
  \end{minipage}
  \begin{minipage}{.32\textwidth}
    \centering
    \includegraphics[width=0.95\textwidth]{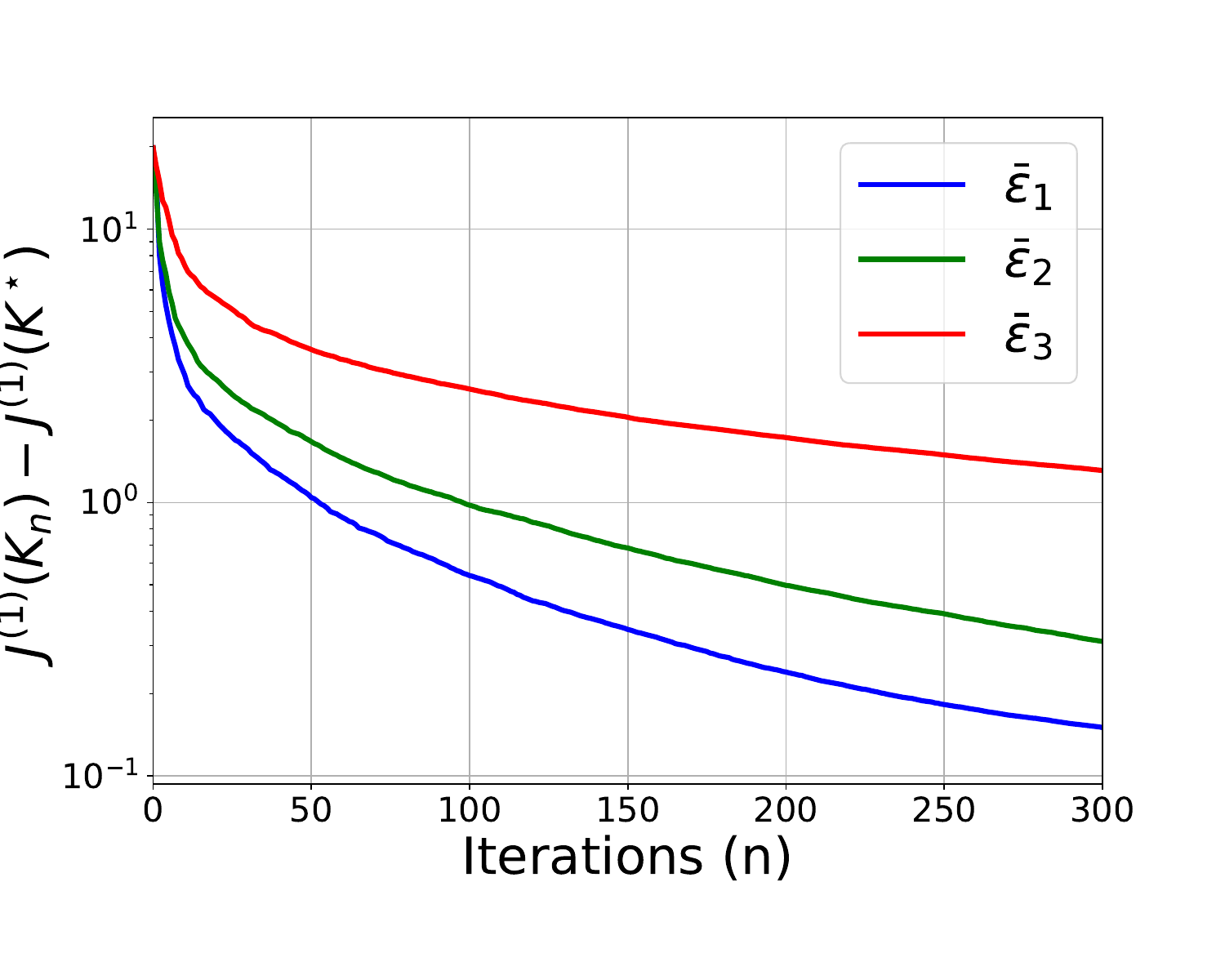}
  \end{minipage}
  \begin{minipage}{.32\textwidth}
    \centering
    \includegraphics[width=0.95\textwidth]{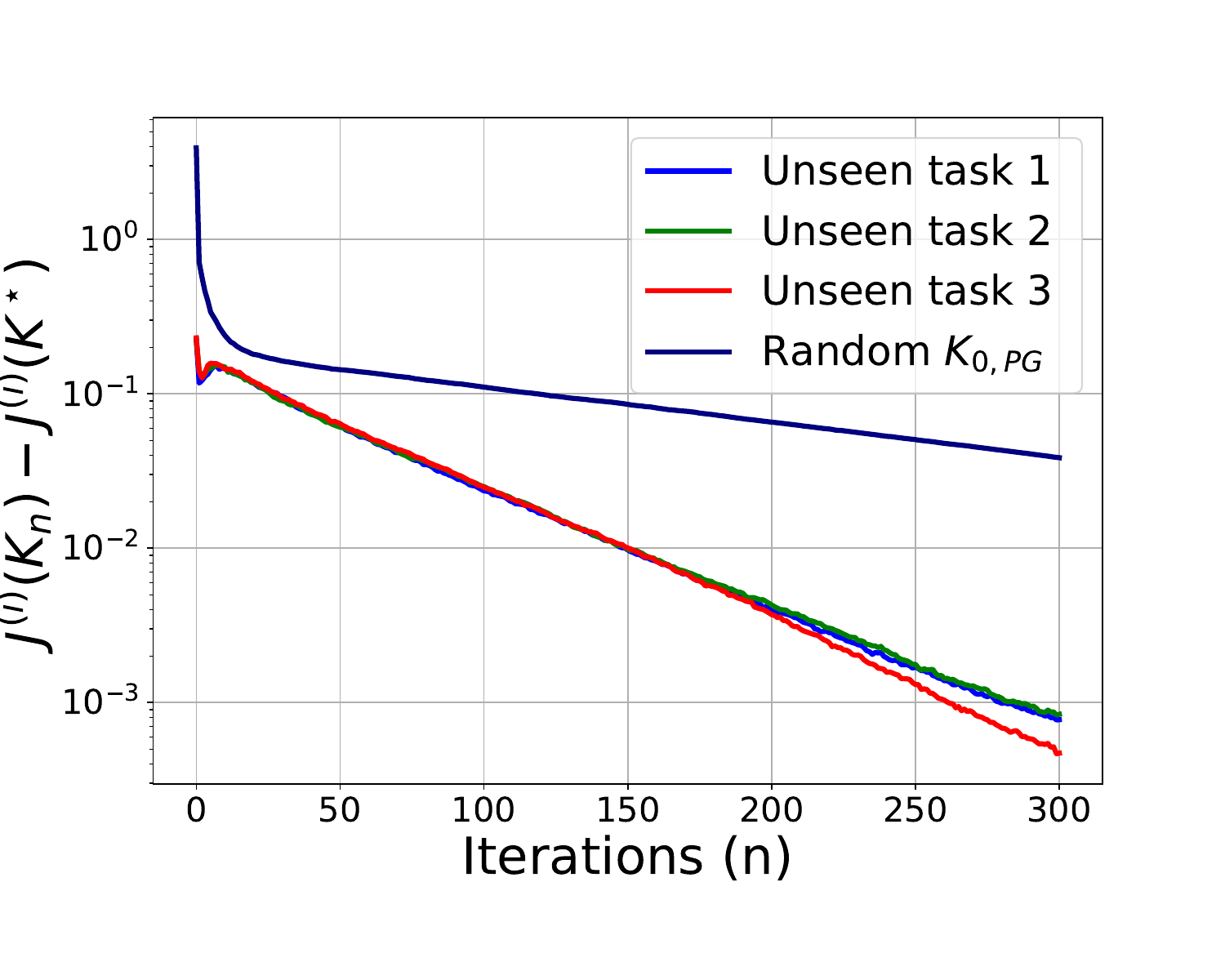}
  \end{minipage}
  \caption{Cost gap between the learned the task-specific optimal controller with respect to iteration. (left) Convergence of the \texttt{MAML-LQR}. (middle) \texttt{MAML-LQR}, $\bar{\epsilon}_1 = (1.2, 1.1, 1.4, 1.2)\times 10^{-3}$, $ \bar{\epsilon}_2 = (1.3, 1.1, 1.4, 1.2)\times 10^{-2}$, $\bar{\epsilon}_3 = (1.7, 1.8, 1.9, 1.7)\times 10^{-2}$. (right) \texttt{PG-LQR} \citep{fazel2018global}.}
  \label{fig:numericals}
\end{figure}

Figure \ref{fig:numericals} (left and middle) depicts the cost gap between the current learned controller and the nominal task (i.e., $\mathcal{T}^{(1)}$) optimal controller with respect to iterations of Algorithm \ref{alg:MAML_model_free}. In alignment with Theorem \ref{theorem:convergence_model_free}, Figure \ref{fig:numericals}-(left) shows that the learned controller closely converges to the nominal task's optimal controller up to a slight bias characterized by $\bar{\epsilon} = (1.2, 1.1, 1.4, 1.2)\times 10^{-3}$. Moreover, Figure \ref{fig:numericals}-(middle) shows that the learned MAML-LQR controller drastically deviates from the nominal task's optimal controller when it faces a significant heterogeneity level. This aligns with Theorem \ref{theorem:convergence_model_free}, which demonstrates that the learned MAML-LQR controller is close to each task-specific optimal controller up to a heterogeneity bias, where the bias increases with $\bar{\epsilon}$.

Figure \ref{fig:numericals}-(right) illustrates the adaptation of the learned MAML-LQR controller to an unseen task drawn from $p(\mathcal{T})$ (i.e., the same task distribution used in the MAML-LQR learning process). With unseen tasks 1, 2 and 3, this figure shows that, by initializing the PG-LQR approach \citep{fazel2018global} from the MAML-LQR learned controller, it takes only a few PG iterations to achieve a controller that is sufficiently close to the unseen tasks' optimal, which is significantly fewer than initializing from a randomly sampled initial stabilizing controller $K_{0,\text{PG}}$. This aligns with Theorem \ref{theorem:convergence_model_free}, showing that the learned controller is close to each task-specific optimal controller.

\vspace{-0.4cm}
\section{Conclusions and Future Work}

We investigated the problem of meta-learning linear quadratic regulators in a heterogeneous and model-free setting, characterizing the stability and convergence of a MAML-LQR approach. We provided theoretical guarantees to ensure task-specific stability under the learned controller for both model-based and model-free settings. We established gradient heterogeneity bounds for three different task heterogeneity cases and offered convergence guarantees showing that the learned controller is close to each task-specific optimal controller up to a task-heterogeneity bias, emphasizing its ability to adapt to unseen tasks. Numerical experiments demonstrated the effect of task heterogeneity on the MAML-LQR approach's convergence and assessed the learned controller's adaptation to unseen tasks. Future work may explore variance-reduced approaches to reduce the variance of the zeroth-order gradient and Hessian estimation to improve the model-free sample complexity further.

\section*{Acknowledgements} 

The authors thank Bruce D. Lee for the valuable comments on this work. Leonardo F. Toso is funded by the Columbia Presidential Fellowship. James Anderson is partially funded by NSF grants ECCS 2144634 and 2231350 and the Columbia Data Science Institute.

\bibliographystyle{abbrvnat}
\bibliography{bibliography}

\newpage

\section{Appendix}

We now provide the technical details/proofs of our main results. For this purpose, we first revisit some auxiliary lemmas and norm inequalities that are instrumental in our stability and convergence analysis of the MAML-LQR approach. Besides the lemmas, we also provide the explicit expressions of the positive polynomials $h_G(K)$, $h_c(K)$, $h_H(K)$, $h_{\Delta}(K)$, $h_{\text {cost}}(K)$, and $h_{\text {grad}}(K)$ that appear in our theoretical guarantees. Lemmas \ref{lemma:Bernstein_gradient} and \ref{lemma:Bernstein_hessian} are used to control the gradient and Hessian estimation errors through a matrix Bernstein-type of inequality (Lemma \ref{lemma:Bernstein}). These lemmas are crucial to provide the stability and convergence guarantees for the model-free setting. Moreover, we provide the experimental setup details utilized to generate the results illustrated in Figure \ref{fig:numericals}.\\

\noindent \textbf{Notation:} Consider a set of $M$ matrices $\{A^{(i)}\}_{i=1}^M$, we denote $\|A\|_{\max} := \max_{i} \|A^{(i)}\|$, and  $\|A\|_{\min} := \min_{i} ||A^{(i)}||$.  The spectral radius of a square matrix is $\rho(\cdot)$.

\subsection{Auxiliary Lemmas and Norm Inequalities} \label{appendix:polynomials}

We first introduce auxiliary lemmas and norm inequalities that are essential in proving our main results. \\

\noindent \textbf{Young's inequality:}  Given any two matrices $A, B \in \mathbb{R}^{n_x\times n_u}$,  for any $\beta>0$, we have

\begin{align}\label{eq:youngs}
\|A+B\|_2^2 &\leq(1+\beta)\|A\|_2^2+\left(1+\frac{1}{\beta}\right)\|B\|_2^2 \leq (1+\beta)\|A\|_F^2+\left(1+\frac{1}{\beta}\right)\|B\|_F^2.
\end{align}

Moreover, given any two matrices $A, B$ of the same dimensions,  for any $\beta>0$, we have
\begin{align}\label{eq:youngs_inner_product}
\langle A, B\rangle & \leq \frac{\beta}{2}\lVert A\rVert_2^2 +\frac{1}{2\beta}\lVert B \rVert_2^2  \leq  \frac{\beta}{2}\lVert A\rVert_F^2 +\frac{1}{2\beta}\lVert B \rVert_F^2.
\end{align}

\noindent \textbf{Jensen's inequality:} Given $M$ matrices $A^{(1)}, \ldots, A^{(M)}$ of identical dimensions, we have that 
\begin{align}\label{eq:sum_expand}
\left\|\sum_{i=1}^M A^{(i)}\right\|_2^2 \leq M \sum_{i=1}^M\left\|A^{(i)}\right\|_2^2, 
\left\|\sum_{i=1}^M A^{(i)}\right\|_F^2 \leq M \sum_{i=1}^M\left\|A^{(i)}\right\|_F^2.
\end{align}

\begin{lemma}(Uniform bounds) Given a LQR task $i\in [M]$ and an stabilizing controller $K \in \mathcal{S}_{\text{ML}}$, the Frobenius norm of gradient $\nabla J^{(i)}(K)$, Hessian $\nabla^2 J^{(i)}(K)$ and control gain $K$ can be bounded as follows:
$$
\|\nabla J^{(i)}(K)\|_F \leq h_G(K), \text{ } \|\nabla^2 J^{(i)}(K)\|_F \leq h_H(K), \text { and } \|K\|_F \leq h_c(K),
$$
with 
\begin{align*}
    h_G(K)&= \frac{J_{\max}(K) \sqrt{\frac{\max_{i}\left\|R^{(i)}+B^{(i)\top} P^{(i)}_K B^{(i)}\right\|\left(J_{\max}(K)-J_{\min}\left(K\right)\right)}{\mu}}}{\min_i\sigma_{\min}(Q^{(i)})},\\
    h_H(K)&= \left(2\|R\|_{\max} + \frac{2\|B\|_{\max}{J}_{\max}(K)}{\mu} + \frac{4\sqrt{2}\tilde{\xi}_{\max}\|B\|_{\max}{J}_{\max}(K)}{\mu}\right)\frac{{J}_{\max}(K)n_u}{\|Q\|_{\min}}, \\
    h_c(K)&=\frac{\sqrt{\frac{\max_{i}\left\|R^{(i)}+B^{(i)\top} P^{(i)}_K B^{(i)}\right\|\left(J_{\max}(K)-J_{\min}\left(K\right)\right)}{\mu}}+\left\|B^{\top} P_K A\right\|_{\max}}{\sigma_{\min}(R)},
\end{align*}
with $\tilde{\xi}_{\max} := \frac{1}{\|Q\|_{\min}}\left( \frac{(1 + \|B\|^2_{\max}){J}_{\max}(K)}{\mu} + \|R\|_{\max} - 1\right)$.
 
\begin{proof} The proof for the uniform bounds of $\|\nabla J^{(i)}(K)\|_F$, $\|K\|_F$ and the gradient domination are provided in \citep{fazel2018global,wang2023model}. In addition, the proof for the uniform bound of $\|\nabla^2 J^{(i)}\|_F$ can be found in \citep[Lemma 7.9]{bu2019lqr}. 
\end{proof}

\end{lemma}

\begin{lemma}(Local smoothness) Given two controller $K, K^{\prime} \in \mathcal{S}_{\text{ML}}$ such that $\|\Delta\| := \|K^{\prime} -K\|  \leq h_{\Delta}(K)\ <\infty$.  The LQR cost, gradient and Hessian satisfy:
\begin{align*}
&\left|J^{(i)}\left(K^{\prime}\right)-J^{(i)}(K)\right| \leq h_{\text {cost}}(K) J^{(i)}(K) \|\Delta\|_F, \\
&\left\|\nabla J^{(i)}\left(K^{\prime}\right)-\nabla J^{(i)}(K)\right\|_F \leq h_{\text {grad}}(K)\|\Delta\|_F, \\
&\left\|\nabla^2 J^{(i)}\left(K^{\prime}\right)-\nabla^2 J^{(i)}(K)\right\|_F \leq h_{\text {hess }}(K)\|\Delta\|_F,
\end{align*}
for any $i \in [M]$, where 
\begin{align*}
    h_{\Delta}(K) &= \frac{\max_i \sigma_{\min}(Q^{(i)}) \mu}{4 ||B||_{\max} J_{\max}(K)\left(\left\|A-BK\right\|_{\max}+1\right)}, \\
    h_{\text {cost}}(K) &= \frac{4 \operatorname{Tr}\left(\Sigma_0\right)J_{\max}(K)\|R\|_{\max}}{\mu\min_i \sigma_{\min}\left(Q^{(i)}\right)}\left(\|K\|+\frac{h_{\Delta}(K)}{2}+\|B\|_{\max}\|K\|^2 \left(\left\|A -BK\right\|_{\max}+1\right) \nu(K)\right), \\
    h_{\text {hess }}(K) &= \sup _{\|X\|_F=1}2(h_1(K)+ 2h_2(K))\|X\|^2_F,\\
    h_{\text {grad}}(K) &= 4\left(\frac{J_{\max}(K)}{\min_i\sigma_{\min}(Q)}\right)\Big[\|R\|_{\max}+\|B\|_{\max}\left(\|A\|_{\max}+\|B\|_{\max}\left(\|K\|+h_{\Delta}(K)\right)\right)\nonumber\\
 &\times\left(\frac{h_{\text {cost }}(K) J_{\max}(K)}{\operatorname{Tr}\left(\Sigma_0\right)}\right)+ \|B\|^2_{\max}\frac{J_
{\max}(K)}{\mu}\Big]\\
&+8\left(\frac{J_{\max}(K)}{\min_i\sigma_{\min}(Q)}\right)^2\left(\frac{\|B\|_{\max}\left(\left\|A-BK\right\|_{\max}+1\right)}{\mu}\right)h_0(K).
\end{align*}
with $\nu(K) = \frac{J_{\max}(K)}{\min_i \sigma_{\min}(Q^{(i)})\mu}$, $h_0(K) = \sqrt{\frac{\max_{i}\left\|R^{(i)}+B^{(i)\top} P^{(i)}_K B^{(i)}\right\|\left(J_{\max}(K)-J_{\min}\left(K\right)\right)}{\mu}},$ and 
\begin{align*}
    h_1(K) &= h_3(K)\|B\|^2_{\max} \frac{J_{\max}(K)}{\min_i\sigma_{\min }(Q^{(i)})} + \Tilde{\mu} h_4(K)\|B\|_{\max} \frac{J_{\max}(K)}{\mu} + h_4(K) \max_i \operatorname{Tr}(R^{(i)}),\\
    h_2(K) &= \|B\|_{\max} J_{\max}(K)\left(\frac{h_6(K) h_4(K) \max_i \operatorname{Tr}\left(A^{(i)} - B^{(i)}K\right)}{\mu} + \|B\|_{\max} h_6(K) \Tilde{\mu} \nu(K)\right. \\ &\left.\hspace{9.7cm}+\frac{\Tilde{\mu} h_7(K)}{\min_i\sigma_{\min }(Q^{(i)})}\right),\\
    h_3(K) &= 6\left(\frac{J_{\max}(K)}{\min_i \sigma_{\min }(Q^{(i)})}\right)^2\|K\|^2\|R\|_{\max}\|B\|_{\max}(\|A-B K\|_{\max}+1 )\\
    &+6\left(\frac{J_{\max}(K)}{\min_i \sigma_{\min }(Q^{(i)})}\right)\|K\|\|R\|_{\max},\\
    h_4(K) &= 4\left(\frac{J_{\max}(K)}{\min_i\sigma_{\min }(Q^{(i)})}\right)^2 \frac{\|B\|_{\max}(\|A-B K\|_{\max}+1)}{\mu},\\
    h_6(K) &= \sqrt{\frac{1}{\min_i \sigma_{\min }(Q^{(i)})}\left(\|R\|_{\max}+\frac{1+\|B\|^2_{\max}}{\mu} J_{\max}(K)\right)-1},\\
    h_7(K) &= 4\left(\nu(K) h_8(K) +8\nu^2(K)\|B\|_{\max}\left(\left\|A - BK\right\|_{\max}+1\right) h_9(K)\right),\\
    h_8(K)&=\|R\|_{\max}+\|B\|^2_{\max} \frac{J_{\max}(K)}{\mu} +\left(\|B\|_{\max}\|A\|_{\max}+\|B\|^2_{\max}\|K\|_{\max}\right) h_3(K),\\
    h_9(K)&= 2\left(\|R\|_{\max}\|K\|+\|B\|_{\max}\|A-B K\|_{\max} \frac{J_{\max}(K)}{\mu}\right).
\end{align*}
and $\Tilde{\mu} = 1+\frac{\mu}{h_{\Delta}(K)}$.

\begin{proof}
The proof of the local smoothness for the cost and gradient are provided in \citep[Appendix F]{wang2023model}, whereas the local smoothness of the Hessian is provided in \citep[Lemma 7]{musavi2023convergence}.    
\end{proof}
\end{lemma}

\begin{lemma} (Matrix Bernstein Inequality) \citep[Lemma B.5]{gravell2020learning} Let $\left\{Z_i\right\}_{i=1}^m$ be a set of $m$ independent random matrices of dimension $d_1 \times d_2$ with $\mathbb{E}\left[Z_i\right]=Z$, $\left\|Z_i-Z\right\| \leq B_r$ almost surely, and maximum variance 

$$\max \left(\left\|\mathbb{E}\left(Z_i Z_i^{\top}\right)-Z Z^{\top}\right\|,\left\|\mathbb{E}\left(Z_i^{\top} Z_i\right)-Z^{\top} Z\right\|\right) \leq \sigma_r^2,$$ 
and sample average $\widehat{Z}:=\frac{1}{m} \sum_{i=1}^m Z_i$. Let a small tolerance $\epsilon \geq 0$ and small probability $0 \leq \delta \leq 1$ be given. If
$$
m \geq \frac{2 \min \left(d_1, d_2\right)}{\epsilon^2}\left(\sigma_r^2+\frac{B_r \epsilon}{3 \sqrt{\min \left(d_1, d_2\right)}}\right) \log \left[\frac{d_1+d_2}{\delta}\right]
$$
$\text { then } \mathbb{P}\left[\|\widehat{Z}-Z\|_F \leq \epsilon\right] \geq 1-\delta \text {.}$   
\label{lemma:Bernstein}
\end{lemma}

\begin{lemma}\label{lemma:Bernstein_gradient} (Gradient Estimation) Given a small tolerance $\epsilon$ and probability $0 < \delta < 1$ and suppose the smoothing radius and number of samples is selected according to

\begin{align}\label{eq:smoothing_radius}
    r \leq \underline{h}^1_r\left(\frac{\epsilon}{2}\right):=\min \left\{\underline{h}_{\Delta}, \frac{1}{\bar{h}_{\text {cost }}}, \frac{\epsilon}{2 \bar{h}_{\text {grad }}}\right\},
\end{align}

\begin{align}\label{eq:number_of_samples}
m \geq \bar{h}^1_m\left(\frac{\epsilon}{2}, \delta\right) :=\frac{8 \min\{n_x,n_u\}}{\epsilon^2}\left(\sigma_r^2+\frac{B_r \epsilon}{6 \sqrt{\min\{n_x,n_u\}}}\right) \log \left[\frac{n_x+n_u}{\delta}\right],
\end{align}
where 
\begin{align*}
    &B_r := 2n_xn_u \bar{h}_{\text {cost}} \bar{J}_{\max} + \frac{\epsilon}{2} + \bar{h}_G,\;\ \sigma^2_r := \left(2n_xn_u \bar{h}_{\text {cost}} \bar{J}_{\max}\right)^2 + \left(\frac{\epsilon}{2} + \bar{h}_G\right)^2,
\end{align*}
where $\bar{J}_{\max} := \sup_{K \in {\mathcal{S}_{\text{ML}}}} \max_{i} J^{(i)}(K)$. Then, for any task $\mathcal{T}^{(i)}$, $\|\widehat{\nabla} J^{(i)}(K) - \nabla J^{(i)}(K)\|_F \leq \epsilon$, with probability $1-\delta$.\\

\begin{proof} To control $\|\widehat{\nabla} J^{(i)}(K) - {\nabla} J^{(i)}(K)\|$, for any task $\mathcal{T}^{(i)}$, we use triangle inequality to write

\begin{align*}
    \|\widehat{\nabla} J^{(i)}(K) - {\nabla} J^{(i)}(K)\|_F \leq \|\widehat{\nabla} J^{(i)}(K) - {\nabla} J^{(i)}_r(K)\|_F + \|{\nabla} J^{(i)}(K) - {\nabla} J^{(i)}_r(K)\|_F
\end{align*}
where ${\nabla} J^{(i)}_r(K) = \mc{E}_{U \sim \mathbb{S}_r}[\nabla J^{(i)}(K+U)]$, $\widehat{\nabla}J^{(i)}(K) = \frac{n_x n_u}{2r^2m}\sum_{l=1}^{m}(J^{(i)}(K+U_l) - J^{(i)}(K-U_l))U_l$, which implies $\|K+U_l- K \|_F = \|U_l\|_F = r$ and $\|K-U_l-K \|_F = \|U_l\|_F = r$. Therefore, by selecting $r \leq h_{\Delta}(K)$ we have that both $\left|J\left(K^{\prime}\right)-J^{(i)}(K)\right| \leq h_{\text {cost}}(K) J^{(i)}(K) \|\Delta\|_F$ and $\left\|\nabla J^{(i)}\left(K^{\prime}\right)-\nabla J^{(i)}(K)\right\|_F \leq h_{\text {grad}}(K)\|\Delta\|_F$ are verified (Lemma \ref{Lemma:lipschitz}). Then, by enforcing $r \leq \frac{1}{h_{\text {cost}}(K)}$ we have 

\begin{align*}
|J^{(i)}(K+U) - J^{(i)}(K) | &\leq  J^{(i)}(K) \to J^{(i)}(K+U) \leq 2J^{(i)}(K), \\
|J^{(i)}(K-U) - J^{(i)}(K) | &\leq  J^{(i)}(K) \to J^{(i)}(K-U) \leq 2J^{(i)}(K), 
\end{align*}
which ensures the stability of the LQR task under the perturbed controllers. This implies that $J^{(i)}(K+U)$ and $J^{(i)}(K-U)$ are well-defined. To control $\|{\nabla} J^{(i)}(K) - {\nabla} J^{(i)}_r(K)\|_F$ we can select $r\leq \frac{\epsilon}{2h_{\text {grad}}(K)}$ to obtain
\begin{align*}
    \|\nabla J^{(i)}(K+U) - \nabla J^{(i)}(K)\|_F \leq h_{\text {grad}}(K) r \leq \frac{\epsilon}{2},
\end{align*}
thus, we use the fact that $\nabla J^{(i)}_r(K) = \mc{E}_{U\sim \mathbb{S}_r} [\nabla J^{(i)}(K+U)]$ and Jensen's inequality to write 
\begin{align*}
    \|{\nabla} J^{(i)}(K) - {\nabla} J^{(i)}_r(K)\|_F &= \|\mc{E}_{U\sim \mathbb{S}_r} (\nabla J^{(i)}(K+U) - \nabla J^{(i)}(K))\|_F \\
    &\leq \mc{E}_{U\sim \mathbb{S}_r} \|\nabla J^{(i)}(K+U) - \nabla J^{(i)}(K)\|_F \leq \frac{\epsilon}{2}.
\end{align*}

On the other hand, to control $\|\widehat{\nabla} J^{(i)}(K) - {\nabla} J^{(i)}_r(K)\|_F$ we use the Lemma \ref{lemma:Bernstein}. For this purpose, we denote $Z ={\nabla} J^{(i)}_r(K) = \mc{E}_{U \sim \mathbb{S}_r} [\widehat{\nabla} J^{(i)}(K)]$ and each individual sample is denoted by $$Z_i = \frac{n_x n_u}{2r^2} (J^{(i)}(K+U) - J^{(i)}(K-U))U,$$
then we can write 
\begin{align*}
    \|Z_i\|_F = \left\|\frac{n_x n_u}{2r^2} (J(K+U) - J(K-U))U\right\|_F \leq \frac{n_x n_u}{2r}|J^{(i)}(K+U) - J^{(i)}(K-U)|,
\end{align*}
and use Lemma \ref{Lemma:lipschitz} to write 
\begin{align*}
    |J^{(i)}(K+U) - J^{(i)}(K-U)| &\leq h_{\text {cost}}(K) J^{(i)}(K- U) \|2U\| = 2 h_{\text {cost}}(K) J^{(i)}(K- U) r,
\end{align*}
and obtain 
\begin{align*}
    \|Z_i\|_F \leq n_xn_u h_{\text {cost}}(K) J^{(i)}(K).
\end{align*}

\quad Now we can also control $\|Z\|_F$ by using triangle inequality 

\begin{align*}
    \|Z\|_F &= \|{\nabla} J^{(i)}_r(K)\|_F = \|{\nabla} J^{(i)}_r(K) -\nabla J^{(i)}(K) + \nabla J^{(i)}(K)\|_F \leq \|{\nabla} J^{(i)}_r(K) -\nabla J^{(i)}(K)\|_F\\
    &+ \|\nabla J(K)\|_F \leq \frac{\epsilon}{2} + h_G(K),
\end{align*}
and we can use triangle inequality to write 
\begin{align*}
    \|Z_i - Z \|_F \leq \|Z_i\|_F + \|Z\|_F \leq B_r := n_xn_u h_{\text {cost}}(K) J(K) + \frac{\epsilon}{2} + h_G(K),
\end{align*}
and for the variance term we have
\begin{align*}
    \|\mc{E}(Z_iZ_i^\top)  - ZZ^\top\|_F &\leq \|\mc{E}(Z_iZ_i^\top)\|_F + \|ZZ^\top\|_F\\
    & \leq \max_{Z_i} (\|Z_i\|_F)^2 + \|Z\|^2_F\\
    &\leq \sigma^2_r := \left(n_xn_u h_{\text {cost}}(K) J(K)\right)^2 + \left(\frac{\epsilon}{2} + h_G(K)\right)^2.
\end{align*}
\end{proof}
which completes the proof.  
\end{lemma}

\begin{lemma} \label{lemma:Bernstein_hessian} (Hessian Estimation) Given a small tolerance $\epsilon$ and probability $0<\delta < 1$ and suppose the smoothing radius $r$ and number of samples $m$ is selected according to

\begin{align}\label{eq:smoothing_radius_hess}
    r \leq \underline{h}^2_r\left(\frac{\epsilon}{2}\right):=\min \left\{\underline{h}_{\Delta}, \frac{1}{\bar{h}_{\text {cost }}}, \frac{\epsilon}{2 \bar{h}_{\text {hess }}}\right\},
\end{align}

\begin{align}\label{eq:number_of_samples_hess}
m \geq \bar{h}^2_m\left(\frac{\epsilon}{2}, \delta\right) :=\frac{8 n_u}{\epsilon^2}\left(\bar{\sigma}_r^2+\frac{\bar{B}_r \epsilon}{6 \sqrt{n_u}}\right) \log \left[\frac{2n_u}{\delta}\right],
\end{align}
where 
\begin{align*}
    &\bar{B}_r := \frac{n^2_u(n_u+r^2)}{r}\bar{h}_{\text {cost}} \bar{J}_{\max} + \frac{\epsilon}{2} + \bar{h}_H,\\
    &\bar{\sigma}^2_r := \left(\frac{n^2_u(n_u+r^2)}{r}\bar{h}_{\text {cost}} \bar{J}_{\max} \right)^2 + \left(\frac{\epsilon}{2} + \bar{h}_H\right)^2,
\end{align*}
then, for any task $\mathcal{T}^{(i)}$, $\|\widehat{\nabla}^2 J^{(i)}(K) - \nabla^2 J^{(i)}(K)\|_F \leq \epsilon$ with probability $1-\delta$.\\

\begin{proof} To control $\|\widehat{\nabla}^2 J^{(i)}(K) - {\nabla}^2 J^{(i)}(K)\|$ we can use triangle inequality to write

\begin{align*}
    \|\widehat{\nabla}^2 J^{(i)}(K) - {\nabla}^2 J^{(i)}(K)\|_F \leq \|\widehat{\nabla}^2 J^{(i)}(K) - {\nabla}^2 J^{(i)}_r(K)\|_F + \|{\nabla}^2 J^{(i)}(K) - {\nabla}^2 J^{(i)}_r(K)\|_F,
\end{align*}
where $J^{(i)}_r(K) = \mc{E}_{U \sim \mathbb{S}_r}[J^{(i)}(K+U)]$, $\nabla^2 J^{(i)}_r(K) = \mc{E}_{U \sim \mathbb{S}_r}[\nabla^2 J^{(i)}(K+U)]$, and estimation $\widehat{\nabla}^2J^{(i)}(K) = \frac{n^2_u}{r^2m}\sum_{l=1}^{m}(J^{(i)}(K+U_l) - J^{(i)}(K))(U_lU_l^\top - I_{n_u})$. Then, similar to the proof of Lemma \ref{lemma:Bernstein_gradient}, by selecting the smoothing radius according to $r \leq \min\left\{\frac{1}{h_{\text {cost}}(K)},h_{\Delta}(K) \right\}$ we have that $J(K+U)$ is well-defined and the local smoothness of the cost and Hessian are satisfied. Then, to control $\|{\nabla}^2 J^{(i)}(K) - {\nabla}^2 J_r(K)\|_F$ we select $r\leq \frac{\epsilon}{2h_{\text {hess}}(K)}$ to obtain
\begin{align*}
    \|\nabla^2 J^{(i)}(K+U) - \nabla^2 J^{(i)}(K)\|_F \leq h_{\text {hess}}(K) r \leq \frac{\epsilon}{2},
\end{align*}
and by Jensen's inequality we can write 
\begin{align*}
    \|{\nabla}^2 J^{(i)}(K) - {\nabla}^2 J^{(i)}_r(K)\|_F &= \|\mc{E}_{U\sim \mathbb{S}_r} (\nabla^2 J^{(i)}(K+U) - \nabla^2 J^{(i)}(K))\|_F \\
    &\leq \mc{E}_{U\sim \mathbb{S}_r} \|\nabla^2 J^{(i)}(K+U) - \nabla^2 J^{(i)}(K)\|_F \leq \frac{\epsilon}{2}.
\end{align*}

Therefore, by using Lemma \ref{lemma:Bernstein} to control $\|\widehat{\nabla}^2 J^{(i)}(K) - {\nabla}^2 J^{(i)}_r(K)\|_F$, we have that
$Z ={\nabla}^2 J^{(i)}_r(K) = \mc{E}_{U \sim \mathbb{S}_r} [\widehat{\nabla}^2 J^{(i)}(K)]$ where $Z_i = \frac{n_x n_u}{2r^2} (J^{(i)}(K+U) - J^{(i)}(K))(UU^\top - I_{n_u})$, then we have
\begin{align*}
    \|Z_i\|_F = \|\frac{n^2_u}{2r^2} (J^{(i)}(K+U) - J^{(i)}(K))(UU^\top - I_{n_u})\|_F \leq \frac{n^2_u}{2r^2}|J^{(i)}(K+U) - J^{(i)}(K)|(n_u+r^2),
\end{align*}
and use Lemma \ref{Lemma:lipschitz} to write 
\begin{align*}
    |J^{(i)}(K+U) - J^{(i)}(K)| &\leq h_{\text {cost}}(K) J^{(i)}(K) \|U\|_F = h_{\text {cost}}(K) J^{(i)}(K) r,
\end{align*}
and obtain 
\begin{align*}
    \|Z_i\|_F \leq \frac{n^2_u(n_u+r^2)}{r}h_{\text {cost}}(K) J^{(i)}(K).
\end{align*}
and
\begin{align*}
    \|Z\|_F = \|{\nabla}^2 J_r(K)\|_F &= \|{\nabla}^2 J^{(i)}_r(K) -\nabla^2 J^{(i)}(K) + \nabla^2 J^{(i)}(K)\|_F \leq \|{\nabla}^2 J^{(i)}_r(K) -\nabla^2 J^{(i)}(K)\|_F \\
    &+ \|\nabla^2 J^{(i)}(K)\|_F\leq  \frac{\epsilon}{2} + h_H(K),
\end{align*}
which implies
\begin{align*}
    \|Z_i - Z \|_F \leq \|Z_i\|_F + \|Z\|_F \leq \bar{B}_r := \frac{n^2_u(n_u+r^2)}{r}h_{\text {cost}}(K) J^{(i)}(K) + \frac{\epsilon}{2} + h_H(K),
\end{align*}
and 
\begin{align*}
    \|\mc{E}(Z_iZ_i^\top)  - ZZ^\top\|_F &\leq \|\mc{E}(Z_iZ_i^\top)\|_F + \|ZZ^\top\|_F\\
    & \leq \max_{Z_i} (\|Z_i\|_F)^2 + \|Z\|^2_F\\
    &\leq \bar{\sigma}^2_r := \left(\frac{n^2_u(n_u+r^2)}{r}h_{\text {cost}}(K) J^{(i)}(K)\right)^2 + \left(\frac{\epsilon}{2} + h_H(K)\right)^2.
\end{align*}
for the variance term, which completes the proof. 
\end{proof}

\end{lemma}

\subsection{Proof of Lemma \ref{lemma:gradient_heterogeneity}} \label{appendix:proof_gradient_het}

For the simplicity of exposure, we derive the gradient heterogeneity bound for the more general task heterogeneity setting, i.e., system and cost heterogeneity which spans the other cases by setting $\epsilon_1 = \epsilon_2 = 0$ (cost heterogeneity) and $\epsilon_3 = \epsilon_4 = 0$ (system heterogeneity). 

To derive the expression of the gradient heterogeneity bound, we start by first recalling the expression of the gradient of the LQR cost for $\nabla J^{(i)}(K)$ and $\nabla J^{(j)}(K)$,
 \begin{align*}
        \nabla J^{(i)}(K)=2E^{(i)}_K\Sigma^{(i)}_{K},\text{ and } \nabla J^{(j)}(K)=2E^{(j)}_K{\Sigma}^{(j)}_{K},
\end{align*}
where $E^{(i)}_K = R^{(i)}K - B^{\top} P^{(i)}_K(A^{(i)} - B^{(i)} K), \;\ \Sigma^{(i)}_K= \mc{E}_{x^{(i)}_0\sim \mathcal{X}_0}\sum_{t=0}^\infty x_t x^{\top}_t,$ which can be used to write
\begin{align*}
    \| \nabla J^{(i)}(K) - \nabla J^{(j)}(K) \| &=  2\|E^{(i)}_K\Sigma^{(i)}_K - E^{(j)}_K\Sigma^{(j)}_K \|\\
    &\leq 2\left( \|E^{(i)}_K -E^{(j)}_K\| \|\Sigma^{(i)}_K\| + \|E^{(j)}_K\| \|\Sigma^{(i)}_K -\Sigma^{(j)}_K\|\right).
\end{align*}

From \citep[Lemma 13]{fazel2018global} and the definition of $E^{(i)}_K$, we can write
\begin{align*}
    \|\Sigma^{(i)}_K\| \leq \frac{{J}_{\max}(K)}{\min_i \sigma_{\text{min}}(Q^{(i)})} , 
   \quad  \|E^{(i)}_K\|\leq \|R\|_{\max}\|K\| +  \frac{J_{\max}(K)(\|B\|^2_{\max}\|K\| + \|B\|_{\max}\|A\|_{\max}) }{\mu}.
\end{align*}
where $J_{\max}(K) := \max_i J^{(i)}(K)$. Then, let us first control $\|\Sigma^{(i)}_K -\Sigma^{(j)}_K\|$. For this purpose, similar to \citep{wang2023model}, which considers only the system heterogeneity setting, we have
\begin{align*}
    \|\Sigma^{(i)}_K -\Sigma^{(j)}_K\| &= \|\mathcal{T}^{(i)}_K(\Sigma_0)  - \mathcal{T}^{(j)}_K(\Sigma_0) \| \leq 2 \|\mathcal{T}^{(i)}_K\|^2 \|\mathcal{F}^{(i)}_K - \mathcal{F}^{(j)}_K\|\|\Sigma_0\|\\
    &\leq 4\left(\frac{J_{\max}(K)}{\min_i \sigma_{\min}(Q^{(i)})\mu}\right)^2(\epsilon_1 + \epsilon_2\|K\|)\|A - BK\|_{\max} \|\Sigma_0\|\\
    & = 4\epsilon_1\nu^2(K)\|A - BK\|_{\max} \|\Sigma_0\| + 4\epsilon_2\nu^2(K)\|K\|\|A - BK\|_{\max} \|\Sigma_0\|, 
\end{align*}
and 
\begin{align*}
  \mathcal{F}^{(i)}_K(X)&:=(A^{(i)}-B^{(i)} K) X(A^{(i)}-B^{(i)} K)^{\top},\\ \mathcal{T}^{(i)}_K(X) &:=\sum_{t=0}^{\infty}(A^{(i)}-B^{(i)} K)^t X\left[(A^{(i)}-B^{(i)} K)^{\top}\right]^t,  
\end{align*}
being linear operators on some symmetric matrix $X$. We now proceed to control $\|E^{(i)}_K -E^{(j)}_K \|$. To do so, we can first write

\begin{align*}
    \| E^{(i)}_K -E^{(j)}_K \|  \leq \underbrace{\epsilon_4 \|K\|}_{T_1} + \underbrace{\|B^{(i)\top} P^{(i)}_K(A^{(i)} - B^{(i)}K) - B^{(j)\top} P^{(j)}_K(A^{(j)} - B^{(j)}K)\|}_{T_2},
\end{align*}
where $T_1$ comes from the cost heterogeneity, and $T_2$ is similar to the system heterogeneity term in \citep{wang2023model}. We proceed to control $T_2$, by following some similar manipulations as detailed in \citep[Appendix E.1]{wang2023model}. Therefore, we have 
\begin{align*}
    T_2 &\leq \|B\|_{\max}\|P_K\|_{\max}(\epsilon_1 + \epsilon_2\|K\|) + \|B\|_{\max}\|P^{(i)}_K - P_K^{(j)}\|\|A - BK\|_{\max}\\
    &+ \epsilon_2 \|A - BK\|_{\max}\|P_K\|_{\max},
\end{align*}
with $\|P_K\|_{\max} \leq \frac{J_{\max}(K)}{\mu}$. To control $\|P^{(i)}_K - P_K^{(j)}\|$ we can first write,
\begin{align*}
    &\|P^{(i)}_K - P_K^{(j)}\| = \|\mathcal{T}^{(i)}_K(Q^{(i)} + K^\top R^{(i)} K) -  \mathcal{T}^{(j)}_K(Q^{(j)} + K^\top R^{(j)} K)  \|\\
    &\leq \|\mathcal{T}^{(i)}_K(Q^{(i)} + K^\top R^{(i)} K) -  \mathcal{T}^{(j)}_K(Q^{(i)} + K^\top R^{(i)} K) \|\\
    &+ \underbrace{\|\mathcal{T}^{(j)}_K(Q^{(i)} + K^\top R^{(i)} K) -  \mathcal{T}^{(j)}_K(Q^{(j)} + K^\top R^{(j)} K) \|}_{T_3}\\
    &\leq 2\|\mathcal{T}^{(i)}_K\|^2 \|\mathcal{F}^{(i)}_K - \mathcal{F}^{(j)}_K\|\|Q^{(i)} + K^\top R^{(i)}K\| + T_3\\
    &\leq 4\left(\frac{J_{\max}(K)}{\min_i \sigma_{\min}(Q^{(i)})\mu}\right)^2(\epsilon_1 + \epsilon_2\|K\|)\|A - BK\|_{\max}(\|Q\|_{\max} + \|K\|^2\|R\|_{\max}) + T_3.
\end{align*}

\quad To control $T_3$, we can write 
\begin{align*}
    T_3 &= \|\mathcal{T}^{(j)}_K(Q^{(i)} + K^\top R^{(i)} K) -  \mathcal{T}^{(j)}_K(Q^{(j)} + K^\top R^{(j)} K) \|\\
    &\leq \|\mathcal{T}^{(j)}_K(Q^{(i)} + K^\top R^{(i)} K - Q^{(j)} - K^\top R^{(j)} K )\|\\
    &\leq \frac{J_{\max}(K) (\epsilon_3 + \epsilon_4 \|K\|^2)}{\min_i \sigma_{\min}(Q^{(i)}) \mu},
\end{align*}
then we finally have the following expression for the upper bound of $\|P^{(i)}_K - P_K^{(j)}\|$
\begin{align*}
    \|P^{(i)}_K - P_K^{(j)}\| &\leq 4\left(\frac{J_{\max}(K)}{\min_i \sigma_{\min}(Q^{(i)})\mu}\right)^2(\epsilon_1 + \epsilon_2\|K\|)\|A - BK\|_{\max}(\|Q\|_{\max} + \|K\|^2\|R\|_{\max})\\
    & + \frac{J_{\max}(K) (\epsilon_3 + \epsilon_4 \|K\|^2)}{\min_i \sigma_{\min}(Q^{(i)}) \mu}.
\end{align*}

\quad Therefore, by using the above expression in $T_2$ we can derive the upper bound of $\| E^{(i)}_K -E^{(j)}_K \|$ as follows 

\begin{align*}
    &\| E^{(i)}_K -E^{(j)}_K \| \leq \epsilon_1\|B\|_{\max}\left(\|P_K\|_{\max} +4\|A - BK\|^2_{\max}\left(\|Q\|_{\max}+ \|K\|^2\|R\|_{\max}\right)\nu^2(K)\right)\\
    &+\epsilon_2 \left(\|A-BK\|_{\max}\|P_K\|_{\max} + \|B\|_{\max}\|P_K\|_{\max}\|K\|\right.\\
    &\left.\hspace{4.9cm}+ 4\|K\|\|A - BK\|^2_{\max}\left(\|Q\|_{\max}+ \|K\|^2\|R\|_{\max}\right)\nu^2(K)  \right)\\
    &+  \epsilon_3\|B\|_{\max}\|A - BK\|_{\max} \nu(K) + \epsilon_4 \left(\|K\| + \|K\|^2\|B\|_{\max}\|A - BK\|_{\max} \nu(K)\right),
\end{align*}
which can then be combined with the upper bound of $\|\Sigma^{(i)}_K -\Sigma^{(j)}_K\|$, $\|\Sigma^{(i)}_K\|$, and $\|E^{(j)}_K\|$ to obtain 
\begin{align*}
   \|\nabla J^{(i)}(K) - \nabla J^{(i)}(K)\| \leq \bar{f}_z(\bar{\epsilon}) := \epsilon_1 {h}^1_{\text{het}}(K) + \epsilon_2 {h}^2_{\text{het}}(K) + \epsilon_3 {h}^3_{\text{het}}(K) + \epsilon_4 {h}^4_{\text{het}}(K),
\end{align*}
with, 
\begin{align*}
    &{h}^1_{\text{het}}(K) = \frac{{J}_{\max}(K)}{\min_i \sigma_{\text{min}}(Q^{(i)})} \|B\|_{\max}\left(\|P_K\|_{\max} +4\|A - BK\|^2_{\max}\left(\|Q\|_{\max}+ \|K\|^2\|R\|_{\max}\right)\nu^2(K)\right)\\
    &+ 4\left(\|R\|_{\max}\|K\| +  \frac{J_{\max}(K)(\|B\|^2_{\max}\|K\| + \|B\|_{\max}\|A\|_{\max}) }{\mu}\right)\nu^2(K)\|A - BK\|_{\max} \|\Sigma_0\|,\\
    &{h}^2_{\text{het}}(K) = \frac{{J}_{\max}(K)}{\min_i \sigma_{\text{min}}(Q^{(i)})}\left(\|A-BK\|_{\max}\|P_K\|_{\max} + \|B\|_{\max}\|P_K\|_{\max}\|K\|\right.\\
    &\left.+ 4\|K\|\|A - BK\|^2_{\max}\left(\|Q\|_{\max}+ \|K\|^2\|R\|_{\max}\right)\nu^2(K)  \right)\\
    &+ 4\left(\|R\|_{\max}\|K\| +  \frac{J_{\max}(K)(\|B\|^2_{\max}\|K\| + \|B\|_{\max}\|A\|_{\max}) }{\mu}\right)\nu^2(K)\|K\|\|A - BK\|_{\max} \|\Sigma_0\|,\\
    &{h}^3_{\text{het}}(K) = \frac{{J}_{\max}(K)}{\min_i \sigma_{\text{min}}(Q^{(i)})}\|B\|_{\max}\|A - BK\|_{\max} \nu(K),\\
    &{h}^4_{\text{het}}(K) = \frac{{J}_{\max}(K)}{\min_i \sigma_{\text{min}}(Q^{(i)})}\left(\|K\| + \|K\|^2\|B\|_{\max}\|A - BK\|_{\max} \nu(K)\right).
\end{align*}

\subsection{Proof of Theorem \ref{theorem:stability_model_based}} \label{appendix:stability_model_based}

\textbf{Objective:} Given an initial stabilizing controller $K_0 \in \mathcal{S}_{\text{ML}}$, our aim is to provide the conditions on the step-sizes $\eta_l$, $\eta$ and heterogeneity $\bar{f}_z(\bar{\epsilon})$, such that $K_{n+1} \in \mathcal{S}_{\text{ML}}$ for all $n \in \{0,1,\ldots,N-1\}$. 

We first observe that, if $K\in {\mathcal{S}}_{\text{ML}}$, one may select $\eta_l$ such that $\bar{K} = K - \eta_l \nabla J^{(i)}(K)$ is also in the MAML stabilizing sub-level set ${\mathcal{S}}_{\text{ML}}$. To verify it, we can use the local smoothness property of the LQR cost (Lemma \ref{Lemma:lipschitz}) to obtain

\begin{align*}
       J^{(i)}(\bar{K}) - J^{(i)}(K) &\leq \langle \nabla J^{(i)}(K), \bar{K} - K \rangle + \frac{{h}_{\text{grad}}(K)}{2}\|\bar{K} - K\|^2_F\\
       &\leq -\eta_l\langle \nabla J^{(i)}(K), \nabla J^{(i)}(K) \rangle + \frac{{h}_{\text{grad}}(K)\eta_l^2}{2}\|\nabla J^{(i)}(K)\|^2_F\\
       &\leq \left(\frac{{h}_{\text{grad}}(K)\eta_l^2}{2} - \eta_l\right) \|\nabla J^{(i)}(K)\|^2_F\\
       & \stackrel{(i)}{\leq} -\frac{\eta_l}{2}\|\nabla J^{(i)}(K)\|^2_F,
\end{align*}
where $(i)$ follows from $\eta_l \leq \frac{1}{{h}_{\text{grad}}(K)}$. We now use the gradient domination property of the LQR cost (Lemma \ref{lemma:gradient_domination}) to write 
   \begin{align*}
       J^{(i)}(\bar{K}) - J^{(i)}(K^\star_i) \leq \left(1 - \frac{\eta_l \lambda_i}{2}\right)\left(J^{(i)}(K) - J^{(i)}(K^\star_i)\right),
   \end{align*}
which implies that $J^{(i)}(\bar{K}) - J^{(i)}(K^\star_i) \leq J^{(i)}(K) - J^{(i)}(K^\star_i)$. Then, since $K \in {\mathcal{S}}_{\text{ML}}$ we conclude that $J^{(i)}(\bar{K}) - J^{(i)}(K^\star_i) \leq J^{(i)}(K_0) - J^{(i)}(K^\star_i)$, which leads to $\bar{K} \in {\mathcal{S}}_{\text{ML}}$. Therefore, given $K \in \mathcal{S}_{\text{ML}}$ and by selecting the step-size according to $\eta_l \leq \frac{1}{{h}_{\text{grad}}(K)}$ we have $\bar{K} = K - \eta_l \nabla J^{(i)}(K) \in {\mathcal{S}}_{\text{ML}}$. We now proceed to show that $K_1 \in \mathcal{S}_{\text{ML}}$. For this purpose, we consider the first iteration of the MAML-LQR algorithm (Algorithm \ref{alg:MAML_LQR}) as the base case. Then, at the first iteration, we have 
\begin{align*}
    K_1 = K_0 -  \frac{\eta}{M}\sum_{j=1}^{M} H^{(j)}(K_0) \nabla J^{(j)}(K_0 - \eta_l\nabla J^{(j)}(K_0)),
\end{align*}
and by using the local smoothness of the LQR cost, we can write
\begin{align*}
    &J^{(i)}(K_{1}) - J^{(i)}(K_0) \leq \langle \nabla J^{(i)}(K_0), K_{1} - K_{0} \rangle + \frac{\bar{h}_{\text{grad}}}{2}\|K_{1} - K_0\|^2_F \\
    &=\langle \nabla J^{(i)}(K_0), -\eta \nabla J_{\text{ML}}(K_0) - \eta \nabla J^{(i)}(K_0) + \eta \nabla J^{(i)}(K_0) \rangle + \frac{\bar{h}_{\text{grad}}\eta^2}{2}\|\nabla J_{\text{ML}}(K_0)\|^2_F\\
    & \leq -\frac{\eta}{2}\|\nabla J^{(i)}(K_0)\|^2_F + \frac{\eta}{2}\|\nabla J_{\text{ML}}(K_0) - \nabla J^{(i)}(K_0)\|^2_F + \frac{\bar{h}_{\text{grad}}\eta^2}{2}\|\nabla J_{\text{ML}}(K_0)\|^2_F,
\end{align*}
and by using Young's inequality \eqref{eq:youngs} along with the selection of the step-size $\eta \leq \frac{1}{4\bar{h}_{\text{grad}}}$ we can write
\begin{align*}
      J^{(i)}(K_{1}) - J^{(i)}(K_0) \leq -\frac{\eta}{4}\|\nabla {J}^{(i)}(K_0)\|^2_F + \frac{3\eta}{4}\|\nabla J_{\text{ML}}(K_0) -\nabla J^{(i)}(K_0) \|^2_F.
\end{align*}

As we show in details in the proof of Theorem \ref{theorem:convergence_model_based}, we control the term $\|\nabla J_{\text{ML}}(K_0) -\nabla J^{(i)}(K_0) \|^2_F$ as follows
\begin{align*}
    \|\nabla J_{\text{ML}}(K_0) -\nabla J^{(i)}(K_0) \|^2_F \leq 24
    n^3_uf^2_z(\bar{\epsilon}) + 2\eta^2_l(12\bar{h}_{\text{grad}}^2n^2_u + h^2_H)\|\nabla J^{(i)}(K_0)\|^2_F,
\end{align*}
where we impose $\eta_l \leq \min\left\{\frac{n_u}{\sqrt{2}\bar{h}_H}, \frac{1}{\sqrt{2}\bar{h}_{\text{grad}}}\right\}$. Then, we obtain
\begin{align*}
    J^{(i)}(K_{1}) - J^{(i)}(K_0) \leq - \frac{\eta}{8}\|\nabla J^{(i)}(K_0)\|^2_F + 18\eta n^3_uf^2_z(\bar{\epsilon}),
\end{align*}
by selecting $\eta_l \leq \frac{1}{\sqrt{12(12\bar{h}^2_{\text{grad}}n^2_u +\bar{h}_{H}^2 )}}$. Then, we use the gradient domination property of the LQR cost to write
\begin{align*}
    J^{(i)}(K_{1}) - J^{(i)}(K^\star_i)  &\leq \left(1 - \frac{\eta \lambda_i}{8}\right)(J^{(i)}(K_0) - J^{(i)}(K^\star_i)) + 18\eta n_uf^2_z(\bar{\epsilon}),
\end{align*}
and by assuming that 
\begin{align*}
    18 n^3_uf^2_z(\bar{\epsilon}) \leq \min_{i \in [M]} \frac{\lambda_i}{16}\left(J^{(i)}(K_0) - J^{(i)}(K^\star_i) \right), 
\end{align*}
we obtain 
\begin{align*}
    J^{(i)}(K_{1}) - J^{(i)}(K^\star_i)  &\leq \left(1 - \frac{\eta \lambda_i}{16}\right)(J^{(i)}(K_0) - J^{(i)}(K^\star_i)),
\end{align*}
which implies 
\begin{align*}
    J^{(i)}(K_{1}) - J^{(i)}(K^\star_i)  &\leq J^{(i)}(K_0) - J^{(i)}(K^\star_i).
\end{align*}

Therefore, we know that under the above selection of step-sizes and heterogeneity we have that $K_1, K_1 - \eta_l \nabla J^{(i)}(K_1) \in \mathcal{S}_{\text{ML}}$, for all $i \in [M]$. The stability analysis is completed by applying an induction step on the base case for all iterations $n \in \{0,1,\ldots,N-1\}$.

\subsection{Proof of Theorem \ref{theorem:convergence_model_based}} \label{appendix:convergence_model_based}

\textbf{Objective:} Given a set of LQR tasks $\mathcal{T}$, where each task-specific optimal controller is $K^\star_i$. We aim to derive the gap between $K_N$ (i.e., the controller obtained Algorithm \ref{alg:MAML_LQR}) and $K^\star_i$. We do so, by deriving the bound for $J^{(i)}(K_N) - J^{(i)}(K^{\star}_i)$, for all $i \in [M]$. Moreover, we also characterize the gap between $K^\star_{\text{ML}}$ and $K^\star_i$ through $J^{(i)}(K^\star_{\text{ML}}) - J^{(i)}(K^{\star}_i)$. With both bounds in hands, we are able to conclude about the local convergence of Algorithm \ref{alg:MAML_LQR}.

\begin{remark} Observe that, in general, we cannot provide any guarantee about the global convergence for the MAML-LQR approach. This is due to the summation in \eqref{eq:MAML_LQR_cost} over heterogeneous LQR tasks, and it prevent us of having a nice MAML-LQR optimization landscape, i.e., a landscape where gradient domination type of properties are satisfied.
\end{remark}

\noindent \textbf{Controlling $J^{(i)}(K_N) - J^{(i)}(K^{\star}_i)$:} We recall that the MAML-LQR update \eqref{eq:MAML_LQR_updating_rule} is given by $K_{n+1} = K_n - \eta \nabla J_{\text{ML}}(K_n)$ with $\nabla  J_{\text{ML}}(K_n) = \frac{1}{M}\sum_{j=1}^{M}H^{(j)}(K_n)\nabla J^{(j)}(K_n - \eta_l \nabla J^{(j)}(K_n))$. Therefore, to proceed we can use the fact that $K_n, K_{n+1}\in {\mathcal{S}}_{\text{ML}}$ along with the local smoothness property of each LQR task's objective to write

\begin{align*}
    &J^{(i)}(K_{n+1}) - J^{(i)}(K_n) \leq \langle \nabla J^{(i)}(K_n), K_{n+1} - K_{n} \rangle + \frac{\bar{h}_{\text{grad}}}{2}\|K_{n+1} - K_n\|^2_F \\
    &=\langle \nabla J^{(i)}(K_n), -\eta \nabla  J_{\text{ML}}(K_n) - \eta \nabla J^{(i)}(K_n) + \eta \nabla J^{(i)}(K_n) \rangle + \frac{\bar{h}_{\text{grad}}\eta^2}{2}\|\nabla  J_{\text{ML}}(K_n)\|^2_F\\
    & \leq -\frac{\eta}{2}\|\nabla J^{(i)}(K_n)\|^2_F + \frac{\eta}{2}\|\nabla  J_{\text{ML}}(K_n) - \nabla J^{(i)}(K_n)\|^2_F + \frac{\bar{h}_{\text{grad}}\eta^2}{2}\|\nabla  J_{\text{ML}}(K_n)\|^2_F,
\end{align*}
where we use Young's inequality \eqref{eq:youngs} in the last inequality. We re-write $\|\nabla  J_{\text{ML}}(K_n)\|^2_F$ as follows:
\begin{align*}
    \|\nabla  J_{\text{ML}}(K_n)\|^2_F &= \|\nabla  J_{\text{ML}}(K_n) -\nabla J^{(i)}(K_n) + \nabla J^{(i)}(K_n) \|^2_F\\
    &\leq 2 \|\nabla  J_{\text{ML}}(K_n) -\nabla J^{(i)}(K_n) \|^2_F + 2 \|\nabla J^{(i)}(K_n) \|^2_F,
\end{align*}
by using Young's inequality \eqref{eq:youngs}. We use the expression of $\|\nabla  J_{\text{ML}}(K_n)\|^2_F$ to obtain
\begin{align*}
    J^{(i)}(K_{n+1}) - J^{(i)}(K_n) &\leq \left(-\frac{\eta}{2} + \bar{h}_{\text{grad}}\eta^2 \right)\|\nabla J^{(i)}(K_n)\|^2_F\\
    &+ \left(\frac{\eta}{2} + \bar{h}_{\text{grad}}\eta^2\right)\|\nabla  J_{\text{ML}}(K_n) -\nabla J^{(i)}(K_n) \|^2_F\\
    &\leq -\frac{\eta}{4}\|\nabla {J}^{(i)}(K_n)\|^2_F + \frac{3\eta}{4}\|\nabla  J_{\text{ML}}(K_n) -\nabla J^{(i)}(K_n) \|^2_F,
\end{align*}
where the last inequality follows from the selection $\eta \leq \frac{1}{4\bar{h}_{\text{grad}}}$. We now proceed to control the gradient term $\|\nabla  J_{\text{ML}}(K_n) -\nabla J^{(i)}(K_n) \|^2_F$. We do it for a general stabilizing policy $K \in {\mathcal{S}}_{\text{ML}}$ below. 
\begin{align*}
    \|\nabla  J_{\text{ML}}(K) -\nabla J^{(i)}(K) \|^2_F &= \left\|\frac{1}{M}\sum_{j=1}^{M} H^{(j)}(K) \nabla J^{(j)}(\bar{K}) - \nabla J^{(i)}(K)\right\|^2_F\\
    &\leq \frac{1}{M}\sum_{j=1}^M \left\|H^{(j)}(K)\nabla J^{(j)}(\bar{K}) - \nabla J^{(i)}(K)\right\|^2_F,
\end{align*}
with $\bar{K} := K - \eta_l \nabla J^{(j)}(K)$. Then, since $\eta_l \leq \frac{1}{\bar{h}_{\text{grad}}}$ we also have $\bar{K} \in {\mathcal{S}}_{\text{ML}}$. We emphasize that the last inequality is due to Jensen's inequality. We then proceed to bound the term inside the sum.
\begin{align*}
    \|H^{(j)}(K)\nabla J^{(j)}(\bar{K}) - \nabla J^{(i)}(K)\|^2_F &= \|H^{(j)}(K)\nabla J^{(j)}(\bar{K}) - H^{(j)}(K)\nabla J^{(i)}(K)\\
    &\hspace{3.7cm}+ H^{(j)}(K)\nabla J^{(i)}(K) - \nabla J^{(i)}(K)\|^2_F\\
    &\leq 2 \|H^{(j)}(K)\|^2_F\|\nabla J^{(j)}(\bar{K}) - \nabla J^{(i)}(K)\|^2_F\\
    &+ 2\|(H^{(j)}(K) - I)\nabla J^{(i)}(K)\|^2_F\\
    &=2 \|H^{(j)}(K)\|^2_F\|\nabla J^{(j)}(\bar{K}) - \nabla J^{(i)}(K)\|^2_F\\
    &+ 2\eta^2_l\|\nabla^2 J^{(j)}(K)\|^2_F\|\nabla J^{(i)}(K)\|^2_F\\
    &\stackrel{(i)}{\leq} 4(n^2_u + \eta^2_l\bar{h}^2_H)\|\nabla J^{(j)}(\bar{K}) - \nabla J^{(i)}(K)\|^2_F \\
    &+ 2\eta^2_l\bar{h}^2_H \|\nabla J^{(i)}(K)\|^2_F,
\end{align*}
where $(i)$ is due to the fact $\|H^{(j)}(K)\|^2_F = \|I_{n_u} - \eta_l \nabla^2 J^{(j)}(K)\|^2_F \leq 2(n^2_u+\eta_l^2\|\nabla^2 J^{(j)}(K)\|^2_F)$ along with the uniform bound for the Hessian in Lemma \ref{lemma:uniform_bounds}. Then, we now proceed to control $\|\nabla J^{(j)}(\bar{K}) - \nabla J^{(i)}(K)\|^2_F$.
\begin{align*}
    \|\nabla J^{(j)}(\bar{K}) - \nabla J^{(i)}(K)\|^2_F &= \|\nabla J^{(j)}(\bar{K}) - \nabla J^{(j)}(K) + \nabla J^{(j)}(K) - \nabla J^{(i)}(K)\|^2_F\\
    &\leq 2 \underbrace{\|\nabla J^{(j)}(\bar{K}) - \nabla J^{(j)}(K)\|^2_F}_{\text{grad. Lipschitz}} + 2 \underbrace{\|\nabla J^{(j)}(K) - \nabla J^{(i)}(K)\|^2_F}_{\text{grad. heterogeneity}}\\
    &\leq 2 \bar{h}_{\text{grad}}^2 \eta^2_l \|\nabla J^{(j)}(K)\|^2_F + 2n_uf^2_z(\bar{\epsilon}),
\end{align*}
where in the last term we use the relationship between Frobenius and spectral norm, since the gradient heterogeneity bound was previously derived in the spectral norm. We can also use the fact that $\|\nabla J^{(j)}(K)\|^2_F \leq 2\|\nabla J^{(j)}(K) - \nabla J^{(i)}(K)\|^2_F + 2 \|\nabla J^{(i)}(K)\|^2_F$ to obtain
\begin{align*}
    \|\nabla J^{(j)}(\bar{K}) - \nabla J^{(i)}(K)\|^2_F \leq 2n_uf^2_z(\bar{\epsilon})\left( 2\bar{h}_{\text{grad}}^2\eta^2_l + 1 \right) + 4\bar{h}_{\text{grad}}^2\eta^2_l\|\nabla J^{(i)}(K)\|^2_F,
\end{align*}
which implies 
\begin{align*}
    \|H^{(j)}(K)\nabla J^{(j)}(\bar{K}) - \nabla J^{(i)}(K)\|^2_F &\leq  8(n^2_u + \eta^2_l\bar{h}^2_H)\left( 2\bar{h}_{\text{grad}}^2\eta^2_l + 1 \right)n_uf^2_z(\bar{\epsilon})\\
    &+ \left[16(n^2_u + \eta^2_l\bar{h}^2_H)\bar{h}_{\text{grad}}^2\eta^2_l + 2\eta^2_l\bar{h}^2_H \right] \|\nabla J^{(i)}(K)\|^2_F,\\
    &\stackrel{(i)}{\leq} 24n^3_uf^2_z(\bar{\epsilon}) + 2\eta^2_l(12\bar{h}_{\text{grad}}^2n^2_u + \bar{h}^2_H)\|\nabla J^{(i)}(K)\|^2_F,
\end{align*}
where $(i)$ follows from the selection $\eta_l \leq \min\left\{\frac{n_u}{\sqrt{2}\bar{h}_H}, \frac{1}{\sqrt{2}\bar{h}_{\text{grad}}}\right\} $. Then, we obtain 
\begin{align*}
    \|\nabla  J_{\text{ML}}(K) -\nabla J^{(i)}(K) \|^2_F \leq 24n^3_uf^2_z(\bar{\epsilon}) + 2\eta^2_l(12\bar{h}_{\text{grad}}^2n^2_u + \bar{h}^2_H)\|\nabla J^{(i)}(K)\|^2_F.
\end{align*}
which yields to
\begin{align*}
    J^{(i)}(K_{n+1}) - J^{(i)}(K_n) &\leq \left(-\frac{\eta}{4} + \frac{3\eta\eta^2_l}{2}(12\bar{h}_{\text{grad}}^2n^2_u + {h}_{H}^2) \right)\|\nabla J^{(i)}(K_n)\|^2_F + 18\eta n^3_uf^2_z(\bar{\epsilon})\\
    &\stackrel{(i)}{\leq} - \frac{\eta}{8}\|\nabla J^{(i)}(K_n)\|^2_F + 18\eta n^3_uf^2_z(\bar{\epsilon}),
\end{align*}
where $(i)$ follows from the selection $\eta_l \leq \frac{1}{\sqrt{12(12\bar{h}^2_{\text{grad}}n^2_u +\bar{h}_{H}^2 )}}$. Therefore, we use the gradient domination property of the LQR task cost (Lemma \ref{lemma:gradient_domination}) to write 
\begin{align*}
    J^{(i)}(K_{n+1}) - J^{(i)}(K^\star_i)  &\leq \left(1 - \frac{\eta \lambda_i}{8}\right)(J^{(i)}(K_{n}) - J^{(i)}(K^\star_i)) + 18\eta n^3_uf^2_z(\bar{\epsilon}),
\end{align*}
and by telescoping the above expression over $n =\{0,1, \ldots, N-1\}$ we  finally obtain 
\begin{align*}
    J^{(i)}(K_{N}) - J^{(i)}(K^\star_i)  &\leq \left(1 - \frac{\eta \lambda_i}{8}\right)^N(J^{(i)}(K_0) - J^{(i)}(K^\star_i)) + 18\eta n^3_uf^2_z(\bar{\epsilon})\\
    &=\left(1 - \frac{\eta \lambda_i}{8}\right)^N\Delta^{(i)}_0 + \frac{144}{\lambda_i} n^3_uf^2_z(\bar{\epsilon}),
\end{align*}
and by selecting the number of iterations $N \geq \frac{8}{\eta \lambda_i}\log\left(\frac{\Delta^{(i)}_0}{\epsilon^\prime}\right)$ for a small tolerance $\epsilon^\prime$, we have
\begin{align*}
    J^{(i)}(K_{N}) - J^{(i)}(K^\star_i)  &\leq \epsilon^\prime + \frac{144}{\lambda_i} n^3_uf^2_z(\bar{\epsilon}).
\end{align*}

\noindent \textbf{Controlling} $J^{(i)}(K^\star_{\text{ML}}) - J^{(i)}(K^{\star}_i)$: We first exploit the  the gradient domination property as in \citep[Lemma 11]{fazel2018global} to write 
\begin{align*}
&J^{(i)}(K^\star_{\text{ML}}) - J^{(i)}\left(K^\star_i\right) 
\leq \operatorname{Tr}\left(\Sigma^{(i)}_{K^\star_i} {E_{K^\star_{\text{ML}}}^{(i)\top}}\left(R^{(i)}+{B^{(i)\top}} P^{(i)}_{K^\star_{\text{ML}}} B^{(i)}\right)^{-1} E_{K^\star_{\text{ML}}}^{(i)}\right) \notag\\
& \leq \left\|\Sigma^{(i)}_{K^\star_i}\right\| \operatorname{Tr}\left({E_{K^\star_{\text{ML}}}^{(i)\top}}\left(R^{(i)}+{B^{(i)\top}} P^{(i)}_{K^\star_{\text{ML}}} B^{(i)}\right)^{-1} E_{K^\star_{\text{ML}}}^{(i)}\right) \notag\\
&\leq \frac{\bar{J}_{\max}}{\min_i\sigma_{\min}(Q^{(i)})}\left\|\left(R^{(i)}+B^{(i)\top} P^{(i)}_{K^\star_{\text{ML}}} B^{(i)}\right)^{-1}\right\| \operatorname{Tr}\left({E_{K^\star_{\text{ML}}}^{(i)\top}} E_{K^\star_{\text{ML}}}^{(i)}\right) \notag\\
& \leq \frac{\bar{J}_{\max}}{\min_i \sigma_{\min }(R^{(i)}) \min_i\sigma_{\min}(Q^{(i)})}\operatorname{Tr}\left({E_{K^\star_{\text{ML}}}^{(i)\top}} E_{K^\star_{\text{ML}}}^{(i)}\right) \notag\\
& = \frac{\bar{J}_{\max}}{\min_i\sigma_{\min}(R^{(i)})\min_i\sigma_{\min }(Q^{(i)})}\operatorname{Tr}\left((\Sigma^{(i)}_{K^\star_{\text{ML}}})^{-1} \nabla {J^{(i)\top}(K^\star_{\text{ML}})} \nabla J^{(i)}(K^\star_{\text{ML}}) (\Sigma_{K^\star_{\text{ML}}}^{(i)})^{-1}\right) \notag\\
& \leq  \frac{\bar{J}_{\max}}{\min_i\sigma_{\min }(R^{(i)})\min_i\sigma_{\min }(Q^{(i)}) \sigma_{\min}\left(\Sigma^{(i)}_{K^\star}\right)^2}\operatorname{Tr}\left(\nabla J^{(i)\top}(K^\star_{\text{ML}}) \nabla J^{(i)}(K^\star_{\text{ML}})\right) \notag\\
&\leq \frac{\bar{J}_{\max}}{\mu^2 \min_i\sigma_{\min }(R^{(i)})\min_i\sigma_{\min }(Q^{(i)})} \lVert \nabla J^{(i)} (K^\star_{\text{ML}})\rVert_{F}^2\notag,
\end{align*}
now we proceed to control $\lVert \nabla J^{(i)} (K^\star_{\text{ML}})\rVert_{F}^2$. To do so, we first write
\begin{align*}
    \lVert \nabla J^{(i)} (K^\star_{\text{ML}})\rVert_{F}^2 &\leq \underbrace{2\Big\lVert \nabla J_{\text{ML}}(K^\star_{\text{ML}}) \Big\rVert_{F}^2 }_{= 0} + 2\Big\lVert \nabla J^{(i)}(K^\star_{\text{ML}}) -\nabla J_{\text{ML}}(K^\star_{\text{ML}}) \Big\rVert_{F}^2\\
    &\leq 48n^3_uf^2_z(\bar{\epsilon}) + 4\eta^2_l(12\bar{h}_{\text{grad}}^2n^2_u + h^2_H)\|\nabla J^{(i)}(K^\star_{\text{ML}})\|^2_F,
\end{align*}
where the last inequality follows from the upper bound of $\|\nabla J_{\text{ML}}(K) -\nabla J^{(i)}(K) \|^2_F$. Therefore, we can select $\eta_l\leq   \frac{1}{\sqrt{8(12\bar{h}^2_{\text{grad}}n^2_u +\bar{h}_{H}^2 )}}$ to write 
\begin{align*}
    \lVert \nabla J^{(i)} (K^\star_{\text{ML}})\rVert_{F}^2
    \leq 96n^3_uf^2_z(\bar{\epsilon}),
\end{align*}
which implies
\begin{align*}
&J^{(i)}(K^\star_{\text{ML}})-J^{(i)}\left(K^\star_i\right) 
\leq \frac{96\bar{J}_{\max}n^3_uf^2_z(\bar{\epsilon})}{\mu^2 \min_i\sigma_{\min }(R^{(i)})\min_i\sigma_{\min }(Q^{(i)})}, 
\end{align*}
\noindent and completes the proof. 


\subsection{Proof of Theorem \ref{theorem:stability_model_free}} \label{appendix:stability_model_free}

We now move our attention to the model-free setting, where both gradient and Hessian are estimated through a zeroth-order estimation with two-point estimation. Similar to the model-based setting, we begin the stability analysis by setting the base case as the first iteration Algorithm \ref{alg:MAML_model_free}. Before proceeding to show that $K_1 \in \mathcal{S}_{\text{ML}}$, we can first prove that given an initial stabilizing controller $ K_0 \in \mathcal{S}_{\text{ML}}$, one may select $\eta_l$, $r$ and $m$ to ensure that $\widehat{K}_0 := K_0 - \eta_l\widehat{\nabla} J^{(i)}(K_0) \in {\mathcal{S}_{\text{ML}}}$. To do so, we use the local smoothness of the LQR cost (Lemma \ref{Lemma:lipschitz}) to write 

\begin{align*}
     J^{(i)}(\widehat{K}_0) &- J^{(i)}(K_0) \leq \langle \nabla J^{(i)}(K_0), \widehat{K}_0 - K_{0} \rangle + \frac{\bar{h}_{\text{grad}}}{2}\|\widehat{K}_0 - K_0\|^2_F \\
     &= \langle \nabla J^{(i)}(K_0), - \eta_l\widehat{\nabla} J^{(i)}(K_0) \rangle + \frac{\bar{h}_{\text{grad}}\eta^2_l}{2}\|\widehat{\nabla} J^{(i)}(K_0)\|^2_F\\
     &\leq -\frac{\eta_l}{2}\|\nabla J^{(i)}(K_0)\|^2_F + \frac{\eta_l}{2}\|\widehat{\nabla} J^{(i)}(K_0) - \nabla J^{(i)}(K_0)\|^2_F + \frac{\bar{h}_{\text{grad}}\eta^2_l}{2}\|\widehat{\nabla} J^{(i)}(K_0)\|^2_F\\
     &\leq \left( \bar{h}_{\text{grad}}\eta^2_l -\frac{\eta_l}{2} \right)\|\nabla J^{(i)}(K_0)\|^2_F + \left( \bar{h}_{\text{grad}}\eta^2_l + \frac{\eta_l}{2} \right)\|\widehat{\nabla} J^{(i)}(K_0) - \nabla J^{(i)}(K_0)\|^2_F \\
     &\stackrel{(i)}{\leq} -\frac{\eta_l}{4} \|\nabla J^{(i)}(K_0)\|^2_F + \frac{3\eta_l}{4}\|\widehat{\nabla} J^{(i)}(K_0) - \nabla J^{(i)}(K_0)\|^2_F, 
\end{align*}
where $(i)$ follows from the selection $\eta_l \leq \frac{1}{4\bar{h}_{\text{grad}}}$. Therefore, by using the gradient domination property of the LQR cost (Lemma \ref{lemma:gradient_domination}), we obtain 
\begin{align*}
    J^{(i)}(\widehat{K}_0) - J^{(i)}(K^\star_i) &\leq  \left(1 -\frac{\eta_l\lambda_i}{4}\right) \left(J^{(i)}(K_0) - J^{(i)}(K^\star_i) \right) + \frac{3\eta_l}{4}\|\widehat{\nabla} J^{(i)}(K_0) - \nabla J^{(i)}(K_0)\|^2_F \\
    &= \left(1 -\frac{\eta_l\lambda_i}{4}\right) \Delta^{(i)}_0 + \frac{3\eta_l}{4}\|\widehat{\nabla} J^{(i)}(K_0) - \nabla J^{(i)}(K_0)\|^2_F,
\end{align*}
and since $K_0 \in \mathcal{S}_{\text{ML}}$, we can select the smoothing radius $r$ and the number of samples $m$ according to Lemma \ref{lemma:Bernstein_gradient} such that  $\|\widehat{\nabla} J^{(i)}(K_0) - \nabla J^{(i)}(K_0)\|^2_F \leq \psi := \frac{\lambda_i\Delta^{(i)}_0}{6}$. Then, we obtain 
\begin{align*}
    &r \leq \min \left\{\underline{h}^1_r\left(\frac{\sqrt{\psi}}{2}\right) \right\}, \;\ m \geq \max\left\{\bar{h}^1_m\left(\frac{\sqrt{\psi}}{2}, \delta\right) \right\}.
\end{align*}
with probability $1-\delta$. This implies that 
\begin{align*}
    J^{(i)}(\widehat{K}_0) - J^{(i)}(K^\star_i) &\leq  \left(1 -\frac{\eta_l\lambda_i}{8}\right) \left(J^{(i)}(K_0) - J^{(i)}(K^\star_i) \right)\\
    &\leq J^{(i)}(K_0) - J^{(i)}(K^\star_i),
\end{align*}
which yields to $K_0 - \eta_l\widehat{\nabla} J^{(i)}(K_0) \in \mathcal{S}_{\text{ML}}$. We then proceed to show that $K_1 \in \mathcal{S}_{\text{ML}}$. By following similar derivations as in the convergence analysis, we can write

\begin{align*}
    J^{(i)}(K_{1}) - J^{(i)}(K_0) &\leq -\frac{\eta}{8}\|\nabla J^{(i)}(K_0)\|^2_F + 30\eta n^3_uf^2_z(\bar{\epsilon})+ \frac{9\eta}{8}\|\widehat{\nabla} J_{\text{ML}}(K_0) - \nabla J_{\text{ML}}(K_0)\|^2_F,
\end{align*}
by selecting $\eta_l \leq \frac{1}{\sqrt{20(12\bar{h}^2_{\text{grad}}n^2_u +\bar{h}_{H}^2 )}}$ and $\eta \leq \frac{1}{8\bar{h}_{\text{grad}}}$. Therefore, by using the gradient domination property of the LQR cost (Lemma \ref{lemma:gradient_domination}) we can write
\begin{align*}
    J^{(i)}(K_{1}) - J^{(i)}(K^\star_i) &\leq \left(1-\frac{\eta\lambda_i}{8}\right)\left(J^{(i)}(K_0) - J^{(i)}(K^\star_i)\right) + 30\eta n^3_uf^2_z(\bar{\epsilon})\\
    &+ \frac{9\eta}{8}\|\widehat{\nabla} J_{\text{ML}}(K_0) - \nabla J_{\text{ML}}(K_0)\|^2_F,
\end{align*}
where we suppose that the heterogeneity satisfy $f^2_z(\bar{\epsilon}) \leq \frac{\lambda_i\Delta^{(i)}_0}{480 n^3_u}$ to obtain
\begin{align*}
    J^{(i)}(K_{1}) - J^{(i)}(K^\star_i) &\leq \left(1-\frac{\eta\lambda_i}{16}\right)\Delta^{(i)}_0 + \frac{9\eta}{8}\|\widehat{\nabla} J_{\text{ML}}(K_0) - \nabla J_{\text{ML}}(K_0)\|^2_F,
\end{align*}
which implies
\begin{align*}
    J^{(i)}(K_{1}) - J^{(i)}(K^\star_i) &\leq \left(1-\frac{\eta\lambda_i}{32}\right)\Delta^{(i)}_0,
\end{align*}
by selecting the smoothing radius $r$ and number of samples $m$ according to 
\begin{align*}
    &r \leq \min \left\{\underline{h}^1_r\left(\frac{\sqrt{\bar{\psi}}}{2}\right), \underline{h}^2_r\left(\frac{\sqrt{\bar{\psi}}}{2}\right) \right\}, \;\ m \geq \max\left\{\bar{h}^1_m\left(\frac{\sqrt{\bar{\psi}}}{2}, \delta\right),\bar{h}^2_m\left(\frac{\sqrt{\bar{\psi}}}{2}, \delta\right) \right\}.
\end{align*}
which follows from Lemmas \ref{lemma:Bernstein_gradient} and \ref{lemma:Bernstein_hessian}, with $\bar{\psi} := \frac{\lambda_i \Delta^{(i)}_0}{1296}$. This implies that $K_1 \in \mathcal{S}_{\text{ML}}.$ The stability analysis is completed by applying a induction step for all iterations $n \in \{0,1,\ldots,N-1\}$.

\subsection{Proof of Theorem \ref{theorem:convergence_model_free}} \label{appendix:convergence_model_free}

\noindent \textbf{Controlling $J^{(i)}(K_N) - J^{(i)}(K^{\star}_i)$:} We recall that model-free MAML-LQR update is given by $K_{n+1} = K_n - \eta \widehat{\nabla} J_{\text{ML}}(K_n)$ with $\widehat{\nabla}J_{\text{ML}}(K_n) = \frac{1}{M}\sum_{j=1}^{M}\widehat{H}^{(j)}(K_n)\widehat{\nabla}J^{(j)}(K_n - \eta_l \widehat{\nabla} J^{(j)}(K_n))$.

\begin{align*}
    &J^{(i)}(K_{n+1}) - J^{(i)}(K_n) \leq \langle \nabla J^{(i)}(K_n), K_{n+1} - K_{n} \rangle + \frac{\bar{h}_{\text{grad}}}{2}\|K_{n+1} - K_n\|^2_F \\
    &=\langle \nabla J^{(i)}(K_n), -\eta \widehat{\nabla} J_{\text{ML}}(K_n) - \eta \nabla J^{(i)}(K_n) + \eta \nabla J^{(i)}(K_n) \rangle + \frac{\bar{h}_{\text{grad}}\eta^2}{2}\|\widehat{\nabla} J_{\text{ML}}(K_n)\|^2_F\\
    & \leq -\frac{\eta}{2}\|\nabla J^{(i)}(K_n)\|^2_F + \frac{\eta}{2}\|\widehat{\nabla} J_{\text{ML}}(K_n) - \nabla J^{(i)}(K_n)\|^2_F + \frac{\bar{h}_{\text{grad}}\eta^2}{2}\|\widehat{\nabla} J_{\text{ML}}(K_n)\|^2_F\\
    &\leq -\frac{\eta}{4}\|\nabla J^{(i)}(K_n)\|^2_F + \frac{5\eta}{4}\|{\nabla} J_{\text{ML}}(K_n) - \nabla J^{(i)}(K_n)\|^2_F + \frac{9\eta}{8}\|\widehat{\nabla} J_{\text{ML}}(K_n) - \nabla J_{\text{ML}}(K_n)\|^2_F,
\end{align*}
where the last inequality follows from the selection of the step-size $\eta \leq \frac{1}{8\bar{h}_{\text{grad}}}$ along with Young's inequality \eqref{eq:youngs}. Therefore, in contrast to the model-free setting, we need to control both $\|{\nabla} J_{\text{ML}}(K_n) - \nabla J^{(i)}(K_n)\|^2_F$ and $\|\widehat{\nabla} J_{\text{ML}}(K_n) - \nabla J_{\text{ML}}(K_n)\|^2_F$. From our previous derivations for the model-based setting we know that $\|{\nabla} J_{\text{ML}}(K_n) - \nabla J^{(i)}(K_n)\|^2_F$ can be controlled as follows:

\begin{align*}
    \|\nabla J_{\text{ML}}(K) -\nabla J^{(i)}(K) \|^2_F  \leq 24 n^3_uf^2_z(\bar{\epsilon}) + 2\eta^2_l(12\bar{h}_{\text{grad}}^2n^2_u + \bar{h}^2_H)\|\nabla J^{(i)}(K)\|^2_F,
\end{align*}
which can be used to obtain
\begin{align*}
    J^{(i)}(K_{n+1}) - J^{(i)}(K_n) &\leq -\frac{\eta}{8}\|\nabla J^{(i)}(K_n)\|^2_F + 30\eta n^3_uf^2_z(\bar{\epsilon})+ \frac{9\eta}{8}\|\widehat{\nabla} J_{\text{ML}}(K_n) - \nabla J_{\text{ML}}(K_n)\|^2_F,
\end{align*}
by selecting $\eta_l \leq \frac{1}{\sqrt{20(12\bar{h}^2_{\text{grad}}n^2_u +\bar{h}_{H}^2 )}}$. Now, to control $\|\widehat{\nabla} J_{\text{ML}}(K_n) - \nabla J_{\text{ML}}(K_n)\|^2_F$ we can first write 

\begin{align*}
    \|\widehat{\nabla} J_{\text{ML}}(K_n) - \nabla J_{\text{ML}}(K_n)\|_F &= \left\|\frac{1}{M}\sum_{i=1}^M \widehat{H}^{(i)}(K_n)\widehat{\nabla} J^{(i)}(K_n - \eta_l\widehat{\nabla} J^{(i)}(K_n))\right.\\
    &\left. \hspace{4cm}- {H}^{(i)}(K_n){\nabla} J^{(i)}(K_n - \eta_l{\nabla} J^{(i)}(K_n))\right\|_F\\
    &\stackrel{(i)}{\leq} \frac{1}{M}\sum_{i=1}^M\|\widehat{H}^{(i)}(K_n)\widehat{\nabla} J^{(i)}(\widehat{K}_n) - {H}^{(i)}(K_n){\nabla} J^{(i)}(\bar{K}_n)\|_F,\\
\end{align*}
where we control $(I):=\|\widehat{H}^{(i)}(K_n)\widehat{\nabla} J^{(i)}(\widehat{K}_n) - {H}^{(i)}(K_n){\nabla} J^{(i)}(\bar{K}_n)\|_F$ as follows: 
\begin{align*}
    & (I) = \|\widehat{H}^{(i)}(K_n)\widehat{\nabla} J^{(i)}(\widehat{K}_n) - H^{(i)}(K_n)\widehat{\nabla} J^{(i)}(\widehat{K}_n)
    + H^{(i)}(K_n)\widehat{\nabla} J^{(i)}(\widehat{K}_n) - {H}^{(i)}(K_n){\nabla} J^{(i)}(\bar{K}_n)\|_F\\
    & = \|\left(\widehat{H}^{(i)}(K_n)- H^{(i)}(K_n)\right)\widehat{\nabla} J^{(i)}(\widehat{K}_n)+ H^{(i)}(K_n)\left(\widehat{\nabla} J^{(i)}(\widehat{K}_n) - {\nabla} J^{(i)}(\bar{K}_n)\right)\|_F\\
    &\leq \|\widehat{H}^{(i)}(K_n)- H^{(i)}(K_n)\|_F\|\widehat{\nabla} J^{(i)}(\widehat{K}_n)\|_F +\| H^{(i)}(K_n) \|_F \|\widehat{\nabla} J^{(i)}(\widehat{K}_n) - {\nabla} J^{(i)}(\bar{K}_n)\|_F\\
    &\stackrel{(i)}{\leq} \|\widehat{\nabla}^2J^{(i)}(K_n)- \nabla^2 J^{(i)}(K_n)\|_F\|\widehat{\nabla} J^{(i)}(\widehat{K}_n)\|_F +(n_u + \eta_l \bar{h}_H) \|\widehat{\nabla} J^{(i)}(\widehat{K}_n) - {\nabla} J^{(i)}(\bar{K}_n)\|_F\\
    &\stackrel{(ii)}{\leq} \eta_l \bar{h}_G\|\widehat{\nabla}^2J^{(i)}(K_n)- \nabla^2 J^{(i)}(K_n)\|_F \\
    &+ \left[ (n_u + \eta_l \bar{h}_H) + \eta_l \|\widehat{\nabla}^2J^{(i)}(K_n)- \nabla^2 J^{(i)}(K_n)\|_F \right] \|\widehat{\nabla} J^{(i)}(\widehat{K}_n) - {\nabla} J^{(i)}(\bar{K}_n)\|_F,
\end{align*}
where $(i)$ follows from the definitions of $H^{(i)}(K_n) = I_{n_u} -\eta \nabla^2 J^{(i)}(K_n)$ and $\widehat{H}^{(i)}(K_n) = I_{n_u} -\eta \widehat{\nabla}^2 J^{(i)}(K_n)$, along with $\|H^{(i)}(K_n)\|_F\leq n_u +\eta_l\bar{h}_H$. $(ii)$ is due to $\|\widehat{\nabla} J^{(i)}(\widehat{K}_n)\|_F \leq \|\widehat{\nabla} J^{(i)}(\widehat{K}_n) - \nabla J^{(i)}(\bar{K}_n)\|_F + \|\nabla J^{(i)}(\bar{K}_n)\|_F$ and $\|\nabla J^{(i)}(\bar{K}_n)\|_F \leq \bar{h}_G$.  Therefore, we proceed to bound $\|\widehat{\nabla} J^{(i)}(\widehat{K}_n) - {\nabla} J^{(i)}(\bar{K}_n)\|_F$. For this purpose, we can write
\begin{align*}
    \|\widehat{\nabla} J^{(i)}(\widehat{K}_n) - {\nabla} J^{(i)}(\bar{K}_n)\|_F &= \|\widehat{\nabla} J^{(i)}(\widehat{K}_n) - \nabla J^{(i)}(\widehat{K}_n) +\nabla J^{(i)}(\widehat{K}_n) - {\nabla} J^{(i)}(\bar{K}_n)\|_F\\
    &\leq \|\widehat{\nabla} J^{(i)}(\widehat{K}_n) - \nabla J^{(i)}(\widehat{K}_n) \|_F + \|\nabla J^{(i)}(\widehat{K}_n) - {\nabla} J^{(i)}(\bar{K}_n)\|_F\\
    &\stackrel{(i)}{\leq} \|\widehat{\nabla} J^{(i)}(\widehat{K}_n) - \nabla J^{(i)}(\widehat{K}_n) \|_F + \bar{h}_{\text{grad}}\eta_l\|\widehat{\nabla} J^{(i)}(K_n) - {\nabla} J^{(i)}(K_n)\|_F,
\end{align*}
where $(i)$ follow from the local smoothness property of the LQR cost (Lemma \ref{Lemma:lipschitz}). Hence, since $K_n, \widehat{K}_n \in {\mathcal{S}_{\text{ML}}}$, we have with $1-\delta$,
\begin{align*}
    \|\widehat{\nabla} J^{(i)}(\widehat{K}_n) - {\nabla} J^{(i)}(\bar{K}_n)\|_F \leq  \epsilon(1+\bar{h}_{\text{grad}}\eta_l),
\end{align*}
with $r$ and $m$ satisfying the conditions in Lemma \ref{lemma:Bernstein_gradient}. 
\begin{align*}
    (I) &\leq \eta_l \bar{h}_G\|\widehat{\nabla}^2J^{(i)}(K_n)- \nabla^2 J^{(i)}(K_n)\|_F\\
    &+ \epsilon \left[ (n_u + \eta_l \bar{h}_H)+ \eta_l \|\widehat{\nabla}^2J^{(i)}(K_n)- \nabla^2 J^{(i)}(K_n)\|_F \right](1+\bar{h}_{\text{grad}}\eta_l),
\end{align*}
and by setting $r$ and $m$ according to the conditions in Lemma \ref{lemma:Bernstein_hessian}, we obtain 
\begin{align*}
    (I) \leq \eta_l \bar{h}_G\epsilon + \epsilon \left[ (n_u + \eta_l \bar{h}_H) + \eta_l \epsilon \right](1+\bar{h}_{\text{grad}}\eta_l),
\end{align*}
which implies that 
\begin{align*}
    \|\widehat{\nabla} J_{\text{ML}}(K_n) - \nabla J_{\text{ML}}(K_n)\|_F &\leq \eta_l \bar{h}_G\epsilon + \epsilon \left[ (n_u + \eta_l \bar{h}_H) + \eta_l \epsilon \right](1+\bar{h}_{\text{grad}}\eta_l) \stackrel{(i)}{\leq} 6n_u\epsilon,
\end{align*}
where $(i)$ follows from the selection $\eta_l \leq \min \left\{ \frac{1}{\bar{h}_G}, \frac{n_u}{\bar{h}_H}, \frac{1}{\bar{h}_{\text{grad}}}, \frac{1}{2}\right\}$ and from the fact that $0< \epsilon < 1$. Therefore, we obtain 

\begin{align*}
    J^{(i)}(K_{n+1}) - J^{(i)}(K_n) &\leq -\frac{\eta}{8}\|\nabla J^{(i)}(K_n)\|^2_F + 30\eta n^3_uf^2_z(\bar{\epsilon})+ \frac{324n_u\eta\epsilon}{8},
\end{align*}
where we can use the gradient domination property of the LQR cost to write 
\begin{align*}
    J^{(i)}(K_{n+1}) - J^{(i)}(K^\star_i) &\leq \left(1-\frac{\eta\lambda_i}{8}\right)\left(J^{(i)}(K_n) - J^{(i)}(K^\star_i)\right) + 30\eta n^3_uf^2_z(\bar{\epsilon})+ \frac{324n_u\eta\epsilon}{8},
\end{align*}
and we can unroll the above expression over $n = \{0,1,\ldots,N-1\}$ to obtain 
\begin{align*}
    J^{(i)}(K_{N}) - J^{(i)}(K^\star_i) &\leq \left(1-\frac{\eta\lambda_i}{8}\right)^N\left(J^{(i)}(K_0) - J^{(i)}(K^\star_i)\right) + \frac{240}{\lambda_i} n^3_uf^2_z(\bar{\epsilon})+ \frac{324n_u\epsilon}{\lambda_i},
\end{align*}
where by selecting $N \geq \frac{8}{\eta \lambda_i}\log\left(\frac{2\Delta^{(i)}_0}{\epsilon^\prime}\right)$ and $\epsilon := \frac{\epsilon^\prime \lambda_i}{648n_u}$, we obtain

\begin{align*}
    J^{(i)}(K_{N}) - J^{(i)}(K^\star_i) &\leq \epsilon^\prime + \frac{240n^3_uf^2_z(\bar{\epsilon})}{\lambda_i}.
\end{align*}

\noindent \textbf{Controlling} $J^{(i)}(K^\star_{\text{ML}}) - J^{(i)}(K^{\star}_i)$: By following similar derivation as in the model-based setting, we obtain 
\begin{align*}
    J^{(i)}(K^\star_{\text{ML}}) - J^{(i)}(K_i^\star) \leq \frac{96\bar{J}_{\max}n^3_uf^2_z(\bar{\epsilon})}{\mu^2 \min_i\sigma_{\min }(R^{(i)})\min_i\sigma_{\min }(Q^{(i)})} 
\end{align*}
which completes the proof.

\section{Experimental Results - Setup Details} \label{appendix:experimental_setup}

We now provide some details on the experimental setup considered in this work.\\

\noindent \textbf{Task Generative Process:} Given a nominal LQR task $(A,B,Q,R)$ we generate $M$ similar but not identical LQR tasks as follows\\

\begin{enumerate}
    \item We first generate random scalar factors $a^{(i)} \sim \mathcal{U}(0,\epsilon_1)$, $b^{(i)} \sim \mathcal{U}(0,\epsilon_2)$, $q^{(i)} \sim \mathcal{U}(0,\epsilon_3)$, and $r^{(i)} \sim \mathcal{U}(0,\epsilon_4)$,  $\forall i \in [M]$, with $\epsilon_1$, $\epsilon_2$,  $\epsilon_3$ and $\epsilon_4$ being the heterogeneity levels.
    \item These random vectors are combined with modification masks $Z_a, Z_q \in \mathbb{R}^{n_x\times n_x}$, $Z_b \in \mathbb{R}^{n_x\times n_u}$, and $Z_r \in \mathbb{R}^{n_u\times n_u}$ to generate multiple LQR tasks $(A^{(i)}, B^{(i)}, Q^{(i)}, R^{(i)})$ for all $i \in [M]$.
    \item The LQR tasks $(A^{(i)}, B^{(i)}, Q^{(i)}, R^{(i)})$ are then constructed according to: $A^{(i)} = A + a^{(i)}Z_a$, $B^{(i)} = B+ b^{(i)}Z_b$, $Q^{(i)} = Q + q^{(i)}Z_q$ and $R^{(i)} = R + r^{(i)}Z_r$. 
\end{enumerate}

\noindent \textbf{Modified Boeing system}  \citep{honglecnotes}: We consider a modification of the Boeing system described in \citep{honglecnotes} to obtain an unstable nominal LQR task described by 

\begin{align*}
   A =  \left[\begin{array}{cccc}
1.22 & 0.03 & -0.02 & -0.32 \\
0.01 & 0.47 & 4.70 & 0.00 \\
0.02 & -0.06 & 0.40 & 0.00 \\
0.01 & -0.04 & 0.72 & 1.55
\end{array}\right], B = \left[\begin{array}{cc}
0.01 & 0.99 \\
-3.44 & 1.66 \\
-0.83 & 0.44 \\
-0.47 & 0.25
\end{array}\right], Q = I_4 , R = I_2. 
\end{align*}

We implement Algorithm \ref{alg:MAML_model_free} with $\mathcal{X}_0 \stackrel{d}{=} \mathcal{N}(0,\frac{1}{4}I_{n_x})$, $\eta_l = \eta = 8\times 10^{-6}$, $r=1\times 10^{-2}$, $m = 20$, and initial stabilizing controller given by

\begin{align*}
    K_0 = \begin{bmatrix}
               0.613 & -1.535 &  0.303 &  0.396 \\
               0.888 &  0.604 & -0.147 & -0.582
    \end{bmatrix}.
\end{align*}

\noindent \textbf{Personalization:} We assess the personalization of Algorithm \ref{alg:MAML_model_free} in Figure \ref{fig:numericals}-(right) by generating $\bar{M}$ tasks for the modified Boeing system and use $80\%$ of them to learn $K_N$ and with the remaining $20\%$ we construct a held-out set of tasks where we randomly sample unseen tasks to assess the generalization of the learned controller.

\end{document}